\documentclass{article}%
\usepackage{amsfonts}
\usepackage{amsmath}
\usepackage{amssymb}
\usepackage{amsxtra}
\usepackage{graphicx}
\usepackage{geometry}
\usepackage{color}
\usepackage{colortbl}%
\providecommand{\U}[1]{\protect\rule{.1in}{.1in}}
\newtheorem{theorem}{Theorem}

\newtheorem{lemma}[theorem]{Lemma}

\newtheorem{remark}[theorem]{Remark}

\numberwithin{equation}{section}
\begin{document}

\author{G. Cardone\\University of Sannio - Department of Engineering \\Corso Garibaldi, 107 - 84100 Benevento, Italy\\email: giuseppe.cardone@unisannio.it
\and S.A. Nazarov\\Institute of Mechanical Engineering Problems\\V.O., Bolshoi pr., 61, 199178, St. Petersburg, Russia.\\email: srgnazarov@yahoo.co.uk
\and C. Perugia\\University of Sannio - DSGA \\via dei Mulini, 59/A - 82100 Benevento, Italy\\email: cperugia@unisannio.it}
\title{A gap in the continuous spectrum of a cylindrical waveguide with a periodic
perturbation of the surface.}
\date{}
\maketitle

\begin{abstract}
It is proved that small periodic singular perturbation of a cylindrical
waveguide surface may open a gap in the continuous spectrum of the Dirichlet
problem for the Laplace operator. If the perturbation period is long and the
caverns in the cylinder are small, the gap certainly opens.

\medskip

 \medskip


\end{abstract}

\section{Spectra of cylindrical and periodic waveguides}

\subsection{The cylindrical waveguide.}

Let $\Omega=\omega\times\mathbb{R}$ be a cylinder with the
cross-section $\omega\subset\mathbb{R}^{2}$ bounded by a simple
closed contour $\partial\omega$ assumed to be $C^{\infty}$ -
smooth for simplicity (cf. Remark \ref{r1.1} below). Interpreting
$\Omega$, e. g., as an acoustic waveguide with the soft wall
$\partial\Omega$, we consider the Dirichlet
problem for the Helmgoltz equation%
\begin{equation}
-\Delta_{x}v\left(  x\right)  =\mu v\left(  x\right)  \text{,
}x\in \Omega\text{, }v\left(  x\right)  =0\text{,
}x\in\partial\Omega\text{,}
\label{1.1}%
\end{equation}
where $\Delta_{x}$ is the Laplacian, $v$ the pressure, and $\mu$ a
spectral parameter, proportional to square of the oscillation frequency.

It is known (cf. \cite{Wh, Wi}) and can be directly verified that,
above a certain cut-off $\mu_{\dagger}>0$, i. e. for
$\mu\geq\mu_{\dagger}$, the
problem $\left(  \ref{1.1}\right)  $ admits a solution in the form%
\begin{equation}
v\left(  x\right)  =\exp\left(  \pm i\zeta z\right)  V\left(
y\right)
\label{1.2}%
\end{equation}
where $i$ is the imaginary unit, $z=x_{3}$ and $y=\left(
y_{1},y_{2}\right)
=\left(  x_{1},x_{2}\right)  $ while%
\begin{equation}
M=\mu-\zeta^{2} \label{1.3}%
\end{equation}
and $V$ are an eigenvalue and the corresponding eigenfunction of
the model
problem in the cross-section%
\begin{equation}
-\Delta_{y}V\left(  x\right)  =MV\left(  y\right)  \text{,
}y\in\omega\text{,
}V\left(  y\right)  =0\text{, }y\in\partial\omega\text{.} \label{1.4}%
\end{equation}
Let $M_{1}$ be the principal eigenvalue in the spectrum of the
problem $\left(  \ref{1.4}\right)  $:
\begin{equation}
0<M_{1}<M_{2}\leq M_{3}\leq...\leq M_{k}...\rightarrow+\infty.\label{1.001}%
\end{equation}
By the maximum principle (see, e.g., \cite{KuGi}), the eigenvalue
$M_{1}$ is simple and the eigenfunction $V_{1}$ can be fixed such
that
\begin{equation}
\left\Vert V_{1};L^{2}\left(  \omega\right)  \right\Vert =1,\text{ }%
V_{1}\left(  y\right)  >0\text{, }y\in\omega,\text{
}\partial_{n}V_{1}\left(
y\right)  <0,\text{ }y\in\partial\omega, \label{1.000}%
\end{equation}
where $\partial_{n}$ stands for differentiation along the outward
normal and
$L^{2}\left(  \omega\right)  $ for the Lebesgue space. If%
\begin{equation}
\mu\geq\mu_{\dagger}=M_{1}, \label{1.5}%
\end{equation}
then $\zeta$ is a real number in $\left(  \ref{1.3}\right)  $,
function $\left(  \ref{1.2}\right)  $ does not grow or vanish as
$z\rightarrow\pm \infty$ and implies a wave which oscillates in
the case $\zeta\neq0$ and stays constant in $z$ for $\zeta=0.$ In
other words, the wave propagation phenomenon occurs above the
cut-off $\mu_{\dagger}$.

The problem $\left(  \ref{1.1}\right)  $ gives rise to the
unbounded positive and self-adjoint operator $A_{\Omega}$ in
$L^{2}\left(  \Omega\right)  $ with
the differential expression $-\Delta_{x}$ and the domain%
\begin{equation}
\mathcal{D}\left(  A_{\Omega}\right)  =H^{2}\left(  \Omega\right)
\cap\mathring{H}^{1}\left(  \Omega;\partial\Omega\right)  . \label{1.6}%
\end{equation}
We use the standard notation for the Sobolev space and the
subspace $\mathring{H}^{1}\left(  \Omega;\partial\Omega\right)  $
of functions in $H^{1}\left(  \Omega\right)  $ satisfying the
Dirichlet conditions in $\left( \ref{1.1}\right)  $.

\begin{figure}
[h]
\begin{center}
\includegraphics[
height=1.2298in, width=3.8164in
]%
{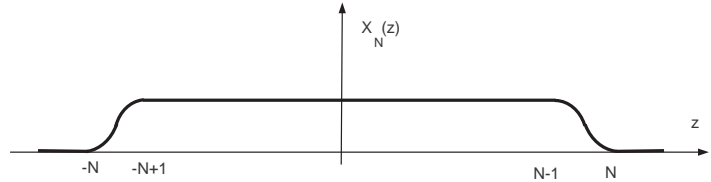}%
\caption{The plateau function}%
\label{fig1}%
\end{center}
\end{figure}

The existence of the nontrivial wave $\left(  \ref{1.2}\right)  $
means that the point $\mu$ belongs to the continuous spectrum
$\sigma_{c}\left( A_{\Omega}\right)  $ of the operator
$A_{\Omega}$. Indeed, multiplying $v$
with the plateau function $X_{N}$ (Figure \ref{fig1}) we see that%
\begin{equation}
\left\Vert X_{N}v;L^{2}\left(  \Omega\right)  \right\Vert
^{2}\geq2\left(
N-1\right)  \text{mes}_{2}\omega,\text{ }\left\Vert \left(  \Delta_{x}%
+\mu\right)  \left(  X_{N}v\right)  ;L^{2}\left(  \Omega\right)
\right\Vert
\leq const,\text{ }N\in\mathbf{%
\mathbb{N}
},\label{Weyl}%
\end{equation}%
where $\mathbf{%
\mathbb{N}
}:=\left\{  1,2,...\right\}  ,$ and, therefore, $\left\{  N^{-\frac{1}{2}%
}X_{h}v\right\}  _{N\in\mathbf{%
\mathbb{N}
}}$ is a singular Weyl sequence of $A_{\Omega}$ at the point $\mu$
whilst $\mu$ belongs to the essential spectrum $\sigma_{e}\left(
A_{\Omega}\right) $ (see, e.g. \cite[\S 9.1]{BiSo}). We emphasize
that, by a general result in
\cite{Ko} (see also \cite[\S 3.1]{NaPl}), the kernel of the mapping%
\begin{equation}
H^{2}\left(  \Omega\right)  \cap\mathring{H}^{1}\left(
\Omega;\partial \Omega\right)  \ni v\mapsto-\left(
\Delta_{x}+\mu\right)  v\in L^{2}\left(
\Omega\right)  \label{1.7}%
\end{equation}
stays finite-dimensional for any $\mu\in\mathbf{%
\mathbb{C}
}$ and, hence, $\sigma_{c}\left(  A_{\Omega}\right)
=\sigma_{e}\left( A_{\Omega}\right)  .$

The set $\mathbf{%
\mathbb{C}
}\backslash\left\{  \mu\in\mathbf{%
\mathbb{C}
:\operatorname{Re}\mu\geq}\mu_{\dagger},\text{ }\operatorname{Im}%
\mu=0\right\}  $ in the complex plane is the resolvent set of the
operator
$A_{\Omega}$ where the inhomogeneous Dirichlet problem%
\begin{equation}
-\Delta_{x}v\left(  x\right)  -\mu v\left(  x\right)  =f\left(
x\right) ,\text{ }x\in\Omega,\text{ }v\left(  x\right)  =0,\text{
}x\in\partial\Omega,
\label{1.8}%
\end{equation}
has a unique solution $v\in H^{2}\left(  \Omega\right)  $ for any
right-hand
side $f\in L^{2}\left(  \Omega\right)  $ and the attendant estimate%
\begin{equation}
\left\Vert v;H^{2}\left(  \Omega\right)  \right\Vert \leq
c_{\mu}\left\Vert
f;L^{2}\left(  \Omega\right)  \right\Vert \label{1.9}%
\end{equation}
ensures that the mapping $\left(  \ref{1.7}\right)  $ is an
isomorphism. In contrast, on the continuous spectrum, mapping
$\left(  \ref{1.7}\right)  $ looses even the Fredholm property
(cf. \cite[Thm. 3.11]{NaPl}) and the inhomogeneous problem $\left(
\ref{1.8}\right)  $ requires a specific formulation involving
radiation conditions at infinity. In the sequel we do
not need such the formulation and refer, e.g. to \cite{Wh, Wi}, \cite[\S 5.3]%
{NaPl} for details.

\subsection{The periodic waveguide, a quasi-cylinder.}

\begin{figure}
[h]
\begin{center}
\includegraphics[
height=0.5094in, width=3.3641in
]%
{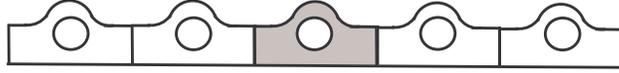}%
\caption{Periodic waveguide}%
\label{fig2}%
\end{center}
\end{figure}

Let $\Pi$ be a domain in $\mathbb{R}^{3}$ with a periodic
cross-section
(Figure \ref{fig2}). More precisely, $\Pi$ is the interior of the union%
\begin{equation}
\overline{\Pi}=%
{\displaystyle\bigcup\nolimits_{j\in\mathbf{\mathbb{Z}}}}
\overline{\varpi_{j}}, \label{1.10}%
\end{equation}%
where $\mathbf{%
\mathbb{Z}
}\mathbf{=}\left\{  0,\pm1,\pm2,...\right\}  ,$
$\varpi_{j}=\left\{  x=\left( y,z\right)  :\left(  y,z-j\right)
\in\varpi\right\}  ,$ and the reference periodicity cell $\varpi$
lies inside the circular cylinder $\left\{
x:\left\vert y\right\vert <R,\text{ }z\in\left(  -\frac{1}{2},\frac{1}%
{2}\right)  \right\}  $ of radius $R>0$ and the unit height. We,
of course, assume that $\Pi$ is a connected set, i.e., a domain.

Due to possible boundary irregularities the Dirichlet problem for
the Helmgoltz equation in the quasi-cylinder $\Pi$ needs the
variational
formulation as the integral identity \cite{Lad}%
\begin{equation}
\left(  \nabla_{x}u,\nabla_{x}v\right)  _{\Pi}=\lambda\left(
u,v\right) _{\Pi},\text{ }v\in\mathring{H}^{1}\left(
\Pi;\partial\Pi\right)  ,
\label{1.11}%
\end{equation}
where $\nabla_{x}$ is the gradient operator, $\left(  \text{
},\text{ }\right)  _{\Pi}$ the natural inner product in
$L^{2}\left(  \Pi\right)  $, and $\lambda$ a spectral parameter.
If the boundary $\partial\Pi$ is smooth, the integral identity
$\left(  \ref{1.11}\right)  $ with arbitrary test function
$v\in\mathring{H}^{1}\left(  \Pi;\partial\Pi\right)  $ is
equivalent to the differential problem of type $\left(
\ref{1.1}\right)  $ in $\Pi$. The
spectral problem reads: To find $\mu\in\mathbf{%
\mathbb{C}
}$ and a nontrivial function $u\in\mathring{H}^{1}\left(  \Omega
;\partial\Omega\right)  $ verifying $\left(  \ref{1.11}\right)  .$

\begin{remark}
\label{r1.1} We could take any open connected and bounded
cross-section $\omega$ of the cylinder $\Omega$ and formulate a
spectral variational problem of type $\left(  \ref{1.11}\right)  $
in $\Omega$. However, the smoothness of $\partial\omega$ will be
used in \S \ 3 for an asymptotic analysis.
\end{remark}

\subsection{The band-gap structure of the essential spectrum in a
quasi-cylinder.}

The left-hand side of $\left(  \ref{1.11}\right)  $ is a positive
continuous form in $\mathring{H}^{1}\left(  \Pi;\partial\Pi\right)
$. According to \cite[\S 10.2]{BiSo}, this form is associated with
a positive self-adjoint unbounded operator $A_{\Pi}$ in
$L^{2}\left(  \Pi\right)  $. If the surface $\partial\Pi$ is
smooth, $A_{\Pi}$ gets the same properties as $A_{\Omega}$ with
the only exception, namely, its essential spectrum\footnote{The authors do not know if it is possible that a segment $\Upsilon_{p}$ collapses into the single point $\Lambda_{p}^+=\Lambda_{p}^-$ which thus becomes an eigenvalue of the operator $A_{\Pi}$. In the case $\Lambda_{p}^+>\Lambda_{p}^-$ for any $p=1,2,...$ the essential spectrum $\sigma_{e}(A_{\Pi})$ coincides with the continuous spectrum $\sigma_{c}(A_{\Pi})$ as in the straight cylinder.} has the
band-gap
structure%
\begin{equation}
\sigma_{e}\left(  A_{\Pi}\right)  =%
{\displaystyle\bigcup\nolimits_{p\in\mathbf{\mathbb{N}}}}
\Upsilon_{p} \label{1.12}%
\end{equation}
where $\Upsilon_{p}$ are closed segments%
\begin{equation}
\Upsilon_{p}=\left[  \Lambda_{p}^{-},\Lambda_{p}^{+}\right]  . \label{1.13}%
\end{equation}
Formulas $\left(  \ref{1.12}\right)  $ and $\left(
\ref{1.13}\right)  $ remain valid without the smoothness
assumption (see \cite{KuchUMN, Kuchbook} and others).

To indicate the segments $\left(  \ref{1.13}\right)  ,$ the model
spectral
problem on the periodicity cell $\varpi$ must be considered%
\begin{equation}
Q_{\eta}\left(  U,V;\varpi\right)  :=\left(  \left(
\nabla_{x}+i\eta e_{3}\right)  U,\left(  \nabla_{x}+i\eta
e_{3}\right)  V\right)  _{\varpi
}=\Lambda\left(  U,V\right)  _{\varpi},\text{ }V\in\mathring{H}_{per}%
^{1}\left(  \varpi;\gamma\right)  , \label{1.14}%
\end{equation}
where $e_{j}$ is the unit vector of the $x_{j}$-axis and $\mathring{H}%
_{per}^{1}\left(  \varpi;\gamma\right)  $ is the subspace of
1-periodic in $z$ functions $V\in H^{1}\left(  \varpi\right)  $
vanishing on the lateral side
$\gamma=\left\{  x\in\partial\varpi:z\in\left(  -\frac{1}{2},\frac{1}%
{2}\right)  \right\}  $ of the cell. If, for certain
$\eta\in\left[
0,2\pi\right)  $ and $\Lambda>0$, the model problem $\left(  \ref{1.14}%
\right)  $ has a nontrivial solution
$U\in\mathring{H}_{per}^{1}\left(
\varpi;\gamma\right)  $, then the Flochet wave%
\begin{equation}
u\left(  y,z\right)  =\exp\left(  i\eta z\right)  U\left(
y,z\right)
\label{1.15}%
\end{equation}
satisfies formally the original problem in the quasi-cylinder
$\Pi$ that is the integral identity $\left(  \ref{1.11}\right)  $
with any test function $v\in C_{c}^{\infty}\left(
\overline{\Pi}\right)  $ (infinitely differentiable functions with
compact supports). One readily constructs from the Flochet wave
the singular Weyl sequence for the operator $A_{\Pi}$ at the point
$\lambda=\Lambda$ with the help of the plateau function drawn in
Figure \ref{fig1} (cf. formulas $\left(  \ref{Weyl}\right)  $).

For any real $\eta$, the sesquilinear form on the left of $\left(
\ref{1.14}\right)  $ is Hermitian, closed and positive. Thus,
problem $\left( \ref{1.14}\right)  $ can be associated with the
unbounded self-adjoint positive operator $\mathcal{A}_{\Pi}\left(
\eta\right)  $ in $L^{2}\left(
\varpi\right)  $ (see again \cite[\S 10.1]{BiSo}). The domain $\mathcal{D}%
\left(  \mathcal{A}_{\Pi}\left(  \eta\right)  \right)  $ is
included into the Sobolev space $H^{1}\left(  \varpi\right)  $
and, therefore, is compactly embedded into $L^{2}\left(
\varpi\right)  $. By \cite[Thm. 10.1.5]{BiSo}, the spectrum of
$\mathcal{A}_{\Pi}\left(  \eta\right)  $ is discrete and forms the
infinitely large sequence
\begin{equation}
0<\Lambda_{1}\left(  \eta\right)  \leq\Lambda_{2}\left(
\eta\right) \leq...\leq\Lambda_{p}\left(  \eta\right)
\leq...\rightarrow+\infty
\label{1.16}%
\end{equation}
of eigenvalues which are listed according to multiplicity. The
functions $\eta\rightarrow\Lambda_{p}\left(  \eta\right)  $ are
continuous (see \cite[Ch. 9]{Kato}) and, by an evident argument,
$2\pi$ - periodic. This means that the endpoints of segments in
$\left(  \ref{1.13}\right)  $ are calculated
as follows:%
\begin{equation}
\Lambda_{p}^{\pm}=\pm\max\left\{  \pm\Lambda_{p}\left(
\eta\right)  :\eta
\in\left[  0,2\pi\right)  \right\}  . \label{1.17}%
\end{equation}

\subsection{The Fourier and Gel'fand transforms.}

Let us comment on the above-mentioned inference. A correspondence
between the problems $\left(  \ref{1.1}\right)  $ in the cylinder
$\Omega$ and $\left( \ref{1.4}\right)  $ in the cross-section
$\omega$ is pointed by the Fourier transform (see \cite{Ko} and
e.g. \cite{NaPl, KoMaRo}). For the quasi-cylinder $\Pi$, one ought
to apply the discrete Fourier transform, namely, the Gel'fand
transform%
\begin{equation}
u\left(  y,z\right)  \longrightarrow\widehat{u}\left(
y,z;\eta\right)
=\frac{1}{\sqrt{2\pi}}%
{\displaystyle\sum\limits_{j\in\mathbf{\mathbb{Z}}}}
\exp\left(  -i\eta\left(  z+j\right)  \right)  u\left(
y,z+j\right)
\label{1.18}%
\end{equation}
(see \cite{Gel} and e.g. \cite{NaPl, Kuchbook}). Note that $\left(
y,z\right)  \in\Pi$ on the left of $\left(  \ref{1.18}\right)  $
but $\left( y,z\right)  \in\varpi$ on the right. The Gel'fand
transform establishes the
isomorphisms%
\[
L^{2}\left(  \Pi\right)  \approx L^{2}\left(  0,2\pi;L^{2}\left(
\varpi\right)  \right)  ,\text{ }H^{l}\left(  \Pi\right)  \approx
L^{2}\left( 0,2\pi;H_{per}^{l}\left(  \varpi\right)  \right)
\]
where $L^{2}\left(  0,2\pi;\mathfrak{B}\right)  $ stands for the
Lebesgue
space of abstract functions,%
\[
\left\Vert U;L^{2}\left(  0,2\pi;\mathfrak{B}\right)  \right\Vert
=\left(
{\displaystyle\int\limits_{0}^{2\pi}}
\left\Vert U\left(  \eta;\mathfrak{B}\right)  \right\Vert
^{2}d\eta\right)
^{\frac{1}{2}}%
\]
and $\mathfrak{B}$ is a Banach space. The corresponding Parceval
theorem
provides the identity%
\[%
{\displaystyle\int\limits_{\Pi}}
u\left(  x\right)  \overline{v\left(  x\right)  }dx=%
{\displaystyle\int\limits_{0}^{2\pi}}
{\displaystyle\int\limits_{\varpi}}
\widehat{u}\left(  x;\eta\right)  \overline{\widehat{v}\left(
x;\eta\right) }dxd\eta,\text{ }u,v\in L^{2}\left(  \Pi\right)  ,
\]
which, together with the formulas%
\[
\widehat{v}\left(  y,0;\eta\right)  =\widehat{v}\left(
y,1;\eta\right) ,\text{ }\partial_{z}\widehat{v}\left(
y,z;\eta\right)  =\widehat {\partial_{z}v}\left(  y,z;\eta\right)
-i\eta\widehat{v}\left( y,z;\eta\right)  ,\text{ }v\in
C_{c}^{\infty}\left(  \overline{\Pi}\right)  ,
\]
indicate the immediate correspondence between the problem $\left(
\ref{1.11}\right)  $ in the quasi-cylinder $\Pi$ and the family
$\left(
\eta\in\left[  0,2\pi\right)  \right)  $ of problems $\left(  \ref{1.14}%
\right)  $ in the periodicity cell $\varpi$.

A result in \cite{na17} (see also \cite{KuchUMN, Kuchbook, NaPl}
and others) demonstrates that the operator of problem $\left(
\ref{1.11}\right)  $ with
the fixed $\lambda\in\mathbf{%
\mathbb{C}
}$ regarding as the mapping
\begin{equation}
\mathring{H}^{1}\left(  \Pi;\partial\Pi\right)
\longrightarrow\mathring
{H}^{1}\left(  \Pi;\partial\Pi\right)  ^{\ast} \label{1.002}%
\end{equation}
is Fredholm if and only if, for any $\eta\in\left[  0,2\pi\right)
$, the problem
\[
\left(  \left(  \nabla_{x}+i\eta e_{3}\right)  U,\left(
\nabla_{x}+i\eta e_{3}\right)  V\right)  _{\varpi}-\lambda\left(
U,V\right)  _{\varpi }=\mathcal{F}\left(  V\right)  ,\text{
}V\in\mathring{H}_{per}^{1}\left( \varpi;\gamma\right)  ,
\]
with the fixed $\lambda\in\mathbf{%
\mathbb{C}
}$ is uniquely solvable with any linear functional
$\mathcal{F\in}\mathring {H}_{per}^{1}\left(  \varpi;\gamma\right)
^{\ast}$ on the space $\mathring {H}_{per}^{1}\left(
\varpi;\gamma\right)$. The
fact mentioned above ensure the segmental structure $\left(  \ref{1.12}%
\right)  $ of the essential spectrum $\sigma_{c}\left(
A_{\Pi}\right)  $ and formulas $\left(  \ref{1.13}\right)  $,
$\left(  \ref{1.17}\right)  $, $\left(  \ref{1.16}\right)  $ for
the segments.

\subsection{Gaps in the essential spectrum.}

The band structure $\left(  \ref{1.12}\right)  $ of the essential
spectrum $\sigma_{e}\left(  A_{\Pi}\right)  $ in the periodic
waveguide $\Pi$ allows for gaps, i.e., intervals on the real
positive semi-axis $\mathbb{R}_{+}$ which lie outside
$\sigma_{e}\left(  A_{\Pi}\right)  $ but have both the endpoints
in $\sigma_{e}\left(  A_{\Pi}\right)  .$ As was commented, such a
gap cannot appear in the essential spectrum $\sigma_e(A_\Omega)=\sigma_c(A_\Omega)$ of the cylindrical waveguide $\Omega$. However, the segments $\left(  \ref{1.13}\right)  $ for the
periodic waveguide can intersect each
other and, as a result, cover the whole ray $\left[  \lambda_{_{\dagger}%
},+\infty\right)  .$ In other words, even a quasi-cylinder can
have the essential spectrum with the only cut-off
$\lambda_{_{\dagger}}$ and no gap.

The main aim of the paper is to show that a small periodic surface
perturbation of the cylinder $\Omega=\omega\times\mathbb{R}$ opens
a gap in the essential spectrum of the corresponding perturbed
quasi-cylinder $\Pi ^{h}$.

In the literature results on opening gaps are mainly related to
periodic media in the whole space $\mathbb{R}^{m}$ of a piece-wise
constant structure described by either scalar differential
equation, or the Maxwell system. We refer to papers \cite{Tat1,
Tat2, Tat3, Tat4, Tat5, Tat6} and reviews \cite{KuchUMN, KuchRew}.
Usually the existence of a gap in the essential spectrum is
established by assuming contrast properties of the media and
selecting or matching the coefficient constants. Results of a different kind are obtained in \cite{japan,SolomyakF,naMaNo} and the present paper, namely coefficients of differential operators are constant and invariable but gaps are opened by varying the shape of the periodicity cell forming a quasicylinder. In \cite{japan,SolomyakF} two-dimensional periodic waveguides of thin width are investigated while spatial waveguides with either regular perturbation of a cylindrical boundary, or a periodic nucleation are studied in \cite{naMaNo}.

\section{Opening a gap in the continuous spectrum in the perturbed periodic
waveguide\label{sect2}}

\subsection{Any cylinder is a periodic set.}

\begin{figure}
[h]
\begin{center}
\includegraphics[
height=1.4434in, width=1.8723in
]%
{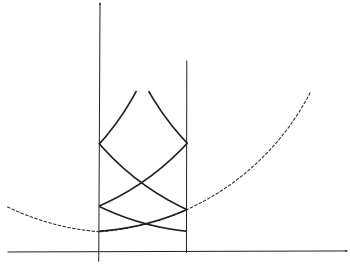}%
\caption{Constructing Floquet waves in the straight cylinder.}%
\label{fig3}%
\end{center}
\end{figure}

\begin{figure}
[h]
\begin{center}
\includegraphics[
height=1.2073in, width=3.4541in
]%
{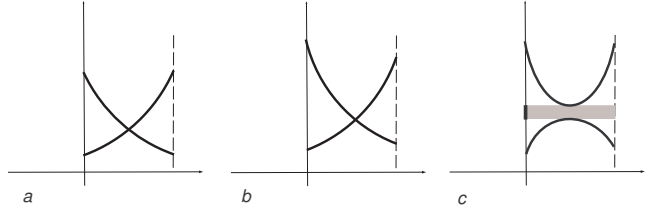}%
\caption{The perturbations (b and c) of the eigenvalue curves (a).}%
\label{fig4}%
\end{center}
\end{figure}

\begin{figure}
[h]
\begin{center}
\includegraphics[
height=1.2592in, width=1.9095in
]%
{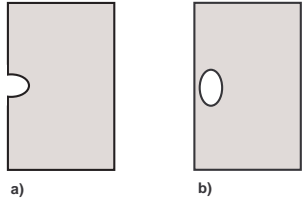}%
\caption{Singular perturbation of the periodicity cell.}%
\label{fig5}%
\end{center}
\end{figure}

The straight cylinder $\Omega=\omega\times\mathbb{R}$ can be
regarded as the quasi-cylinder $\Pi^{0}$ with the periodicity cell
$\varpi^{0}=\omega \times\left(  0,1\right)  .$ The wave $v\left(
x\right)  =\exp\left(  i\zeta z\right)  V\left(  y\right)  $ (see
$\left(  \ref{1.2}\right)  $) turns into
the Floquet wave $\left(  \ref{1.15}\right)  $ with the attributes%
\begin{equation}
\eta=\zeta-2\pi q\left(  \zeta\right)  ,\text{ }U\left(
y,z\right) =\exp\left(  2\pi iq\left(  \zeta\right)  z\right)
V\left(  y\right)  ,
\label{1.19}%
\end{equation}
where $q\left(  \zeta\right)  =\max\left\{  q\in\mathbf{%
\mathbb{Z}
:}2\pi q\leq\zeta\right\}  $. We point out that the factor
$\exp\left(  2\pi iq\left(  \zeta\right)  z\right)  $ is
$2\pi$-periodic in $z$. Thus, each of
the curves%
\begin{equation}
\mu=M_{p}-\zeta^{2}, \label{1.20}%
\end{equation}
forming the continuous spectrum $\sigma_{c}\left(  \Omega\right)
$ (cf. $\left(  \ref{1.3}\right)  $ and $\left(  \ref{1.5}\right)
$), gives rise to
infinite number of pieces%
\begin{equation}
\lambda=M_{p}+\left(  \eta-2\pi q\right)  ^{2},\text{ }q\in\mathbb{%
\mathbb{Z}
}\mathbf{,\ }\eta\in\left[  0,2\pi\right)  \label{1.200}%
\end{equation}
generating segments in $\left(  \ref{1.12}\right)  $ which cover
the whole ray $\left[  M_{1},+\infty\right)  $ as it is shown on
Figure \ref{fig3} for the lowest curve $\left(  \ref{1.20}\right)
$ with $p=1$. In particular, under the assumption
\begin{equation}
M_{1}+\pi^{2}<M_{2} \label{1.21}%
\end{equation}%
(see Remark \ref{r1.2} below) the spectral problem $\left(
\ref{1.14}\right)
$ in the cylindrical cell $\varpi^{0}=\omega\times\left(  -\frac{1}{2}%
,\frac{1}{2}\right)  $ gets the first eigenpairs%
\begin{equation}
\Lambda_{1}^{0}\left(  \eta\right)  =\left\{
\begin{tabular}
[c]{ll}%
$M_{1}+\eta^{2},$ & $\eta\in\left[  0,\pi\right)  ,$\\
$M_{1}+\left(  2\pi-\eta\right)  ^{2},$ & $\eta\in\left(  \pi,2\pi\right)  ,$%
\end{tabular}
\ \ \ \ \ \ \ \ \ \ \right.  \label{1.22}%
\end{equation}%
\begin{equation}
U_{1}^{0}\left(  y,z;\eta\right)  =V_{1}\left(  y\right)  \left\{
\begin{tabular}
[c]{ll}%
$1,$ & $\eta\in\left[  0,\pi\right)  ,$\\
$\exp\left(  -2\pi iz\right)  ,$ & $\eta\in\left(  \pi,2\pi\right)  ,$%
\end{tabular}
\ \ \ \ \ \ \ \ \ \ \right.  \label{1.23}%
\end{equation}
while, at $\eta=\pi$, the eigenvalue%
\begin{equation}
\Lambda_{1}^{0}\left(  \pi\right)  =\Lambda_{2}^{0}\left(
\pi\right)
=M_{1}+\pi^{2} \label{1.24}%
\end{equation}
becomes of multiplicity $2$ and has the eigenfunctions%
\begin{equation}
U_{+}^{0}\left(  y,z\right)  =V_{1}\left(  y\right)  ,\text{
}U_{-}^{0}\left( y,z\right)  =V_{1}\left(  y\right)  \exp\left(
-2\pi iz\right)  .
\label{1.25}%
\end{equation}
It is known (see, e.g., \cite[\S 7.6]{Kato} and
\cite[Ch.9.10]{MaNaPl}) that a small perturbation of the cell
$\varpi^{0}$ prompts perturbations of eigenvalues in $\left(
\ref{1.16}\right)  $. Two situations drawn in Figure \ref{fig4}
may occur for the first couple of eigenvalues and we shall show
that a periodic singular perturbation of the cylinder $\Omega$
(see Figure \ref{fig5}) provides opening a gap in the continuous
spectrum (as indicated by
over-shadowing in Figure \ref{fig4},c).%


\subsection{The singular perturbation of the cylindrical surface.}

\begin{figure}
[h]
\begin{center}
\includegraphics[
height=1.2816in, width=1.4486in
]%
{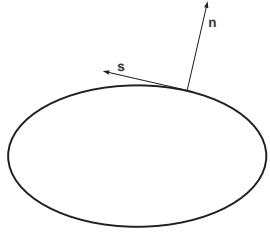}%
\caption{The curvilinear coordinates.}%
\label{fig6}%
\end{center}
\end{figure}

To describe the boundary perturbation of the cylinder
$\Pi^{0}=\Omega$, we introduce in a neighborhood $\Upsilon$ of the
contour $\Gamma=\partial\omega$ the natural curvilinear coordinate
system $\left(  n,s\right)  $ (Figure \ref{fig6}) where $n$ is the
oriented distance to $\Gamma$, $n<0$ inside $\omega\cap\Upsilon$
and $s$ is the arc length on $\Gamma$ evaluated from a point
$\mathcal{O}^{\prime}\in\Gamma$ counter-clockwise so that the
point $\mathcal{O=}\left(  \mathcal{O}^{\prime},0\right)
\in\partial\Pi^{0}$ has
the coordinates $n=0,$ $s=0,$ $z=0$.%

Given a small parameter $h$, we introduce the sets%
\begin{equation}
\theta^{h}=\left\{
x\in\Upsilon\times\mathbb{R}\mathbf{:\xi:=}h^{-1}\left(
n,s,z\right)  \in\theta\right\}  ,\text{
}\varpi^{h}=\varpi^{0}\backslash
\overline{\theta^{h}} \label{1.26}%
\end{equation}
where $\theta$ is a bounded nonempty domain in the half-space $\mathbb{R}%
_{-}^{3}=\left\{  \xi=\left(  \xi_{1},\xi_{2},\xi_{3}\right)  :\xi
_{1}<0\right\}  .$ According to formula $\left(  \ref{1.10}\right)
$, the reference cell $\varpi^{h}$ in $\left(  \ref{1.26}\right)
$ generates the quasi-cylinder $\Pi^{h}$ with a singular
perturbation by the $1$ - periodic family of the caves or
superficial voids $\theta_{j}^{h}=\left\{  x\in
\Upsilon\times\mathbb{R}:h^{-1}\left(  n,s,z-j\right)
\in\theta\right\}  $ (Figure \ref{fig5},a,b).

\begin{remark}
\label{r1.2} We have assumed that the perturbation period $T$ is
equal to $1$. If $T\neq1$, the rescaling $x\longmapsto T^{-1}x$
turns the cylinder $\Omega$ into
$\Omega_{T}=\omega_{T}\times\mathbb{R}$ while the model problem
$\left( \ref{1.4}\right)  $ in the new cross-section
$\omega_{T}=\left\{ y:Ty\in\omega\right\}  $ gets the eigenvalues
$T^{2}M_{k}$ where $M_{k}$ are taken from $\left(
\ref{1.001}\right)  $. Since $M_{1}<M_{2}$, the assumption
$\left(  \ref{1.21}\right)  $ is satisfied in the case%
\begin{equation}
T>\pi\left(  M_{2}-M_{1}\right)  ^{-\frac{1}{2}}. \label{1.27}%
\end{equation}

\end{remark}

\subsection{The boundary layer phenomenon.}

\begin{figure}
[h]
\begin{center}
\includegraphics[
height=1.1554in, width=2.1145in
]%
{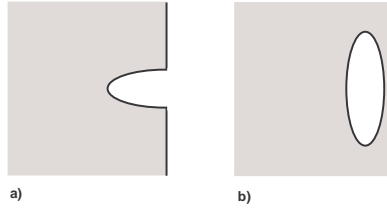}%
\caption{The rescaled perturbation of the boundary.}%
\label{fig7}%
\end{center}
\end{figure}

To examine the behavior of eigenfunctions in the periodicity cell
$\varpi^{h}$ near the boundary perturbation, we need to construct
the boundary layer (see, e.g., \cite{Ilin}, \cite[Ch.
2.9]{MaNaPl}). To this end, we use the stretched coordinates $\xi$
in $\left(  \ref{1.26}\right)  $. Since the Laplacian
$\Delta_{x}$ in the curvilinear coordinates reads%
\begin{equation}
\Delta_{x}=\left(  1+n\varkappa\left(  s\right)  \right)
^{-1}\left( \frac{\partial}{\partial n}\left(  1+n\varkappa\left(
s\right)  \right) \frac{\partial}{\partial
n}+\frac{\partial}{\partial s}\left(  1+n\varkappa \left(
s\right)  \right)  ^{-1}\frac{\partial}{\partial s}\right)
+\frac{\partial^{2}}{\partial z^{2}}, \label{1.28}%
\end{equation}
where $\varkappa\left(  s\right)  $ is the curvature of $\Gamma$
at the point
$s$, we formally have%
\begin{equation}
\Delta_{x}\sim h^{-2}\Delta_{\xi}+h^{-1}\left(  \varkappa\left(
\mathcal{O}^{\prime}\right)
\frac{\partial}{\partial\xi_{1}}-2\varkappa \left(
\mathcal{O}^{\prime}\right)  \xi_{1}\frac{\partial^{2}}{\partial
\xi_{2}^{2}}\right)  +... \label{1.29}%
\end{equation}
Hence, in view of formulae (\ref{1.26}), the coordinate dilation
$x\longmapsto\xi$ leads to the following limit problem%
\begin{equation}
-\Delta_{\xi}w\left(  \xi\right)  =0,\text{ }\xi\in\Theta\text{,
}w\left(
\xi\right)  =g\left(  \xi\right)  ,\text{ }\xi\in\partial\Theta, \label{1.30}%
\end{equation}
in the imperfect half-space (Figure \ref{fig7},a,b)%
\begin{equation}
\Theta=\mathbb{R}_{-}^{3}\backslash\overline{\theta}. \label{1.31}%
\end{equation}%
It is known that, for a sufficiently smooth datum $g$ with a
compact support, problem $\left(  \ref{1.30}\right)  $ has a
unique solution with a finite Dirichlet integral. In the sequel we
need such the decaying solution $W\left( \xi\right)  $ with the
special right-hand side $g\left(  \xi\right)  =-\xi _{1}$ which
vanishes on $\partial\Theta\backslash\partial\theta$ and obeys the
asymptotic form%
\begin{equation}
W\left(  \xi\right)
=-\frac{1}{2\pi}P_{\theta}\frac{\xi_{1}}{\left\vert \xi\right\vert
^{3}}+O\left(  \left\vert \xi\right\vert ^{-3}\right)  ,\text{
}\left\vert \xi\right\vert >R, \label{1.32}%
\end{equation}
where $R>0$ is fixed such that $\left\vert x\right\vert <R$ for
$x\in \overline{\theta}.$

Note that $-\left(  2\pi\left\vert \xi\right\vert ^{3}\right)
^{-1}\xi_{1}$ implies the Poisson kernel and $P_{\theta}>0$ by
virtue of the maximum principle.

\begin{figure}
[h]
\begin{center}
\includegraphics[
height=1.1398in, width=2.1344in
]%
{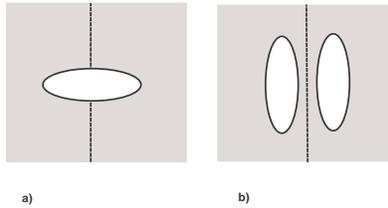}%
\caption{The symmetrization of the boundary perturbation.}%
\label{fig8}%
\end{center}
\end{figure}

\begin{remark}
\label{r1.3} The exterior Dirichlet problem for the symmetrized
set $\theta^{\bullet\bullet}=\left\{  \xi:\left(  -\left\vert
\xi_{1}\right\vert ,\xi_{2},\xi_{3}\right)
\in\overline{\theta}\right\}  $ (cf. Figures \ref{fig7} and
\ref{fig8}) has an intrinsic integral characteristics, the
polarization matrix (see \cite[Appendix G]{PoSe}), which is
extracted from asymptotics at the infinity of the harmonics
$W_{j}$ under the Dirichlet conditions $W_{j}\left(  \xi\right)
=-\xi_{j},$ $\xi\in\partial \theta^{\bullet\bullet}$. The odd
extension of $W$ from $\Theta$ onto
$\mathbb{R}^{3}\backslash\theta^{\bullet\bullet}$ coincides with
$W_{3}$ and, therefore, $P_{\theta}$ is proportional to an entry
in the polarization tensor of $\theta^{\bullet\bullet}$. We call
$P_{\theta}$ the polarization coefficient of the cavity or void
$\theta$ in the half-space.
\end{remark}

%

\subsection{The main result on asymptotics.}

To identify the gap, we need two assertions on eigenvalues and
eigenfunctions
of the auxiliary problem%
\begin{equation}
\left(  \left(  \nabla_{x}+i\eta e_{3}\right)  U^{h},\left(
\nabla_{x}+i\eta e_{3}\right)  V\right)
_{\varpi^{h}}=\Lambda^{h}\left(  U^{h},V\right)
_{\varpi^{h}},\text{ }V\in\mathring{H}_{per}^{1}\left(
\varpi^{h};\gamma
^{h}\right)  , \label{1.36}%
\end{equation}
in the perturbed periodicity cell $\varpi^{h}$ in $\left(
\ref{1.26}\right) $ with the lateral side $\gamma^{h}=\left\{
x\in\partial\varpi^{h}:\left\vert z\right\vert <1/2\right\}  .$ We
enumerate the eigenvalues in the same way as
in $\left(  \ref{1.16}\right)  $:%
\begin{equation}
0<\Lambda_{1}^{h}\left(  \eta\right)  \leq\Lambda_{2}^{h}\left(
\eta\right) \leq...\leq\Lambda_{p}^{h}\left(  \eta\right)
\leq...\rightarrow+\infty.
\label{1.37}%
\end{equation}
However, under the assumption $\left(  \ref{1.21}\right)  $ the
first couple of eigenvalues in $\left(  \ref{1.37}\right)  $ is
denoted by $\Lambda_{\pm }^{h}\left(  \eta\right)  $ while,
according to $\left(  \ref{1.22}\right)  $, we have
$\Lambda_{\pm}^{0}=M_{1}+\left(  \eta-\pi\pm\pi\right)  ^{2}$ in
the limit ($h=0$) problem $\left(  \ref{1.14}\right)  $ in
$\varpi=\varpi^{0}$ and the corresponding eigenfunctions are given
by $\left(  \ref{1.23}\right)  .$

\begin{theorem}
\label{teo1} There exist positive numbers $h_{0},\beta_{0},c_{0}$
such that, for any $h\in\left(  0,h_{0}\right]  $ and $\left\vert
\beta\right\vert \leq\beta_{0}h^{-\frac{5}{4}}$, the first couple
of eigenvalues in $\left( \ref{1.37}\right)  $ of the problem
$\left(  \ref{1.36}\right)  $ on the periodicity cell
$\varpi^{h}$, determined in $\left(  \ref{1.26}\right)  $, takes
the asymptotic form
\begin{equation}
\Lambda_{\pm}^{h}\left(  \pi+\beta h^{3}\right)
=M_{1}+\pi^{2}+h^{3}\left(
\mathcal{P\pm}\sqrt{\mathcal{P}^{2}+4\pi^{2}\beta^{2}}\right)
+\widetilde
{\Lambda}_{\pm}^{h}\left(  \pi+\beta h^{3}\right)  \label{1.33}%
\end{equation}
where the remainder admits the estimate%
\begin{equation}
\left\vert \widetilde{\Lambda}_{\pm}^{h}\left(  \pi+\beta
h^{3}\right)
\right\vert \leq C_{\Lambda}h^{\frac{7}{2}} \label{1.34}%
\end{equation}
and the positive quantity%
\begin{equation}
\mathcal{P=}P_{\theta}\left\vert \partial_{n}V_{1}\left(
\mathcal{O}^{\prime
}\right)  \right\vert ^{2} \label{1.35}%
\end{equation}
is calculated according to $\left(  \ref{1.32}\right)  $ and
$\left( \ref{1.000}\right)  $.
\end{theorem}

The asymptotic formula $\left(  \ref{1.33}\right)  $, $\left(  \ref{1.35}%
\right)  $ will be derived in \S 3 and the remainder estimate
$\left( \ref{1.34}\right)  $ in \S 4. To detect a gap in the
continuous spectrum of the problem $\left(  \ref{1.11}\right)  $
in the quasi-cylinder $\Pi^{h}$, we also prove the following
intelligible inequalities.

\begin{lemma}
\label{lem2} Entries of the eigenvalue sequences $\left(
\ref{1.37}\right)  $ and $\left(  \ref{1.16}\right)  $ of the
auxiliary problems in the cells
$\varpi^{h}$ and $\varpi^{0}$, respectively, are in the relationship%
\begin{equation}
\Lambda_{p}^{0}\left(  \eta\right)  \leq\Lambda_{p}^{h}\left(
\eta\right)
\leq\Lambda_{p}^{0}\left(  \eta\right)  +C_{p}h^{3}, \label{1.39}%
\end{equation}
where $C_{p}$ is independent of $\eta\in\left[  0,2\pi\right)  $
and $h\in\left(  0,h_{0}\right]  .$
\end{lemma}

\textbf{Proof.} Let $\mathcal{A}_{\eta}^{h}$ be a unbounded
operator in $L^{2}\left(  \varpi^{h}\right)  $ generated by the
closed positive Hermitian form $Q_{\eta}\left(
\cdot,\cdot,\varpi^{h}\right)  $ on the left of $\left(
\ref{1.14}\right)  $ (cf. \cite[\S 10.2]{BiSo}). We employ the
max-min
principle (see \cite[Thm. 10.2.2]{BiSo})%
\begin{equation}
\Lambda_{p}^{h}\left(  \eta\right)  =\underset{\mathcal{E}_{p}}{\max}%
\underset{U\in\mathcal{E}_{p}\backslash\left\{  0\right\}  }{\inf}%
\frac{Q_{\eta}\left(  U,U;\varpi^{h}\right)  }{\left\Vert
U;L^{2}\left(
\varpi^{h}\right)  \right\Vert ^{2}},\text{ }p\in\mathbf{%
\mathbb{N}
.} \label{1.40}%
\end{equation}
Here $\mathcal{E}_{p}$ is any subspace in
$\mathring{H}_{per}^{1}\left( \varpi^{h};\gamma^{h}\right)  $ of
co-dimension $p-1$, in particular,
$\mathcal{E}_{1}=\mathring{H}_{per}^{1}\left(
\varpi^{h};\gamma^{h}\right) .$

Let the eigenfunctions $U_{p}^{0}\left(  \cdot,\eta\right)  $
corresponding to $\Lambda_{p}^{0}\left(  \eta\right)  $ satisfy
the normalization and
orthogonality conditions%
\begin{equation}
\left(  U_{p}^{0},U_{q}^{0}\right)  _{\varpi^{0}}=\delta_{p,q},\text{ }%
p,q\in\mathbf{%
\mathbb{N}
},\label{1.42}%
\end{equation}
where $\delta_{p,q}$ stands for Kronecker's symbol. The subspace
$\mathfrak{I}_{p}\subset\mathring{H}_{per}^{1}\left(
\varpi^{h};\gamma
^{h}\right)  $ is spanned over the functions $X_{h}U_{1}^{0},...,X_{h}%
U_{p}^{0}$ while $X_{h}\in C_{per}^{\infty}\left(  \overline{\varpi^{0}%
}\right)  $ is such that
\begin{align}
X_{h}\left(  x\right)   &  =0\text{ \ for \ }\left\vert x-\mathcal{O}%
\right\vert \leq C_{X}h,X_{h}\left(  x\right)  =1\text{ \ for \
}\left\vert
x-\mathcal{O}\right\vert \geq2c_{X}h,\label{1.41}\\
X_{h}\left(  x\right)   &  =0\text{ \ for \ }x\in\theta^{h},\ 0\leq X_{h}%
\leq1,\ \left\vert \nabla_{x}X_{h}\left(  x\right)  \right\vert
\leq ch^{-1}.\ \nonumber
\end{align}
In other words, $X_{h}$ is equal to $1$ everywhere in
$\varpi^{h}$, except in
the vicinity of $\mathcal{O}$, and $X_{h}$ vanishes in the cavern. We have%
\begin{align}
\left(  X_{h}U^{p},X_{h}U^{q}\right)  _{\varpi^{h}} &  =\left(  U^{p}%
,U^{q}\right)  _{\varpi^{0}}+\left(  \left(  1-X_{h}^{2}\right)  U^{p}%
,U^{q}\right)  _{\varpi^{0}}\geq\delta_{p,q}-c_{pq}h^{2}h^{3},\label{1.444}\\
Q_{\eta}\left(  U^{p},U^{q};\varpi^{h}\right)   &  \leq
Q_{\eta}\left( U^{p},U^{q};\varpi^{0}\right)  +c_{pq}\left(
\left\Vert U^{p};H^{1}\left( \Xi_{h}\right)  \right\Vert
+h^{-1}\left\Vert U^{p};L^{2}\left(  \Xi
_{h}\right)  \right\Vert \right)  \cdot\nonumber\\
&  \ \ \ \ \ \ \ \ \ \ \ \ \ \ \ \ \ \ \cdot\left(  \left\Vert U^{q}%
;H^{1}\left(  \Xi_{h}\right)  \right\Vert +h^{-1}\left\Vert
U^{p};L^{2}\left(
\Xi_{h}\right)  \right\Vert \right)  \nonumber\\
&  \leq\Lambda_{p}^{0}\delta_{p,q}+c_{pq}\left(  h^{3}+h^{-2}h^{2}%
h^{3}\right)  .\nonumber
\end{align}
Here come the factors $h^{-1}$ and $h^{2}$ from the
differentiation of $X_{h}$ and the formula
\[
\left\vert U^{p}\left(  x\right)  \right\vert ^{2}\leq
c_{p}\left\vert x-\mathcal{O}\right\vert ^{2}\leq
C_{p}h^{2},\text{ }x\in\Xi_{h},
\]
while $h^{3}$ is order of the volume of the set $\Xi_{h}=$
supp$\left( 1-X_{h}\right)  \supset$ supp$\left\vert
\nabla_{x}X_{h}\right\vert .$

The intersection of the subspaces $\mathcal{E}_{p}$ and
$\mathfrak{I}_{p}$
contains the nontrivial linear combination%
\begin{equation}
\mathcal{U}_{p}=X_{h}%
{\displaystyle\sum\limits_{j=1}^{p}}
a_{p}U^{p},%
{\displaystyle\sum\limits_{j=1}^{p}}
\left\vert a_{p}\right\vert ^{2}=1. \label{1.43}%
\end{equation}
Hence, according to $\left(  \ref{1.40}\right)  $ and $\left(  \ref{1.42}%
\right)  $, we derive that%
\[
\Lambda_{p}^{h}\left(  \eta\right)  \leq\underset{\mathcal{E}_{p}}{\max}%
\frac{Q_{\eta}\left(  \mathcal{U}_{p};\varpi^{h}\right)
}{\left\Vert
\mathcal{U}_{p};L^{2}\left(  \varpi^{h}\right)  \right\Vert ^{2}}\leq\frac{%
{\displaystyle\sum\limits_{j=1}^{p}}
\Lambda_{j}^{0}\left(  \eta\right)  +C_{p}h^{3}}{1-C_{p}h^{3}}%
\]
and the right inequality in $\left(  \ref{1.39}\right)  $ is
proved.

The left inequality can be easily derived by applying the max-min
principle to the operator $\mathcal{A}_{\eta}^{0}$ and extending
eigenfunctions $U_{p}^{h}$ by zero from $\varpi^{h}$ onto
$\varpi^{0}.$

\subsection{Detecting the gap.}

By Lemma \ref{lem2} and formula $\left(  \ref{1.200}\right)  $, we
conclude
that%
\begin{align}
\Lambda_{1}^{h}\left(  \eta\right)   &  \leq M_{1}+\min\left\{
\eta
^{2},\left(  2\pi-\eta\right)  ^{2}\right\}  +C_{1}h^{3},\label{1.44}\\
\Lambda_{2}^{h}\left(  \eta\right)   &  \geq\min\left\{
M_{1}+\max\left\{ \eta^{2},\left(  2\pi-\eta\right)  ^{2}\right\}
,\ M_{2}+\min\left\{ \eta^{2},\left(  2\pi-\eta\right)
^{2}\right\}  \right\}  .\nonumber
\end{align}
Hence, in view of the assumption $\left(  \ref{1.21}\right)  $ we
can choose $\eta_{0}>\left(  2\pi\right)  ^{-1}\max\left\{
C_{1},3\mathcal{P}\right\}  $ and $h_{0}>0$ such that, for
$h\in\left(  0,h_{0}\right]  $ and $\eta
\in\left[  0,\pi-\eta_{0}h^{3}\right)  \cup\left(  \pi+\eta_{0}h^{3}%
,2\pi\right)  $, the interval%
\[
\left(  M_{1}+\pi^{2},M_{1}+\pi^{2}+2\mathcal{P}h^{3}\right)
\]
is free of eigenvalues $\left(  \ref{1.37}\right)  $. In the case
$\left\vert \eta-\pi\right\vert \leq\eta_{0}h^{3}$ we apply
Theorem \ref{teo1} to observe that also the interval
\begin{equation}
\left(  M_{1}+\pi^{2}+C_{\Lambda}h^{\frac{7}{2}},M_{1}+\pi^{2}+2\mathcal{P}%
h^{3}-C_{\Lambda}h^{\frac{7}{2}}\right)  \label{1.45}%
\end{equation}
does not contain the eigenvalues. We emphasize that the endpoints
in $\left( \ref{1.45}\right)  $ are established by the following
inequalities taken from
Theorem \ref{teo1}:%
\begin{align*}
\Lambda_{1}^{h}\left(  \eta\right)   &  \leq
M_{1}+\pi^{2}h^{3}\left(
\mathcal{P-P}\sqrt{1+4\pi^{2}\mathcal{P}^{-2}\eta_{0}^{2}}\right)
+C_{\Lambda}h^{\frac{7}{2}}\leq M_{1}+\pi^{2}+C_{\Lambda}h^{\frac{7}{2}},\\
\Lambda_{2}^{h}\left(  \eta\right)   &  \geq
M_{1}+\pi^{2}h^{3}\left(
\mathcal{P+P}\sqrt{1+4\pi^{2}\mathcal{P}^{-2}\eta_{0}^{2}}\right)
-C_{\Lambda}h^{\frac{7}{2}}\geq
M_{1}+\pi^{2}+2\mathcal{P}h^{3}-C_{\Lambda }h^{\frac{7}{2}}.
\end{align*}
These two facts provide the main result in the paper.

\begin{theorem}
\label{teo3} Under the assumption $\left(  \ref{1.21}\right)  $,
there exist positive numbers $h_{0}$ and $c_{0}$ such that, for
$h\in\left( 0,h_{0}\right]  $, the essential spectrum $\left(
\ref{1.12}\right)  $ of the problem $\left(  \ref{1.11}\right)  $
in the periodic waveguide $\Pi^{h}$ with the periodicity cell
$\varpi^{h}$ in $\left(  \ref{1.26}\right)  $ has a
gap of length $l\left(  h\right)  $,%
\begin{equation}
\left\vert l\left(  h\right)  -2\mathcal{P}h^{3}\right\vert \leq c_{0}%
h^{\frac{7}{2}}, \label{1.46}%
\end{equation}
situated just after the first segment $\Upsilon_{1}^{h}$ in
$\left( \ref{1.12}\right)  $. Here $\mathcal{P}$ is the positive
quantity $\left( \ref{1.35}\right)  $.
\end{theorem}

\begin{figure}
[h]
\begin{center}
\includegraphics[
height=1.5039in, width=1.5887in
]%
{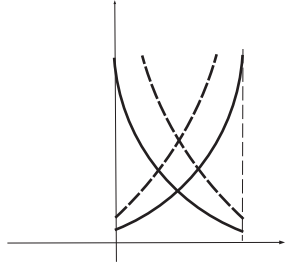}%
\caption{Overlapping of the bands.}%
\label{fig10}%
\end{center}
\end{figure}

We finally mention that under the assumption $M_{1}+\pi^{2}>M_{2}$
opposite to $\left(  \ref{1.21}\right)  ,$ the first and second
segment $\Upsilon_{1}^{h}$ and $\Upsilon_{2}^{h}$ intersect (cf.
Figure \ref{fig10} where two dotted curves correspond to $\left(
\ref{1.200}\right)  $ with $p=2$ and $q=0,1$) and therefore, the
gap discovered in Theorem \ref{teo3} does not occur. This
conclusion readily follows from the rough estimate $\left(
\ref{1.39}\right) $ for the perturbated eigenvalue
$M_{2}+\min\left\{  \eta^{2},\left(
2\pi-\eta\right)  ^{2}\right\}  $ in the auxiliary problem in $\varpi^{h}$.%

If $M_{2}=M_{1}+\pi^{2}$, then the gap is still open because $\Lambda_{2}%
^{h}\left(  0\right)  >M_{2}$ and Theorem \ref{teo1} is still
valid. At the same time, the gap length is $O\left(  h^{3}\right)
$ only in the case when $M_{2}$ is simple and
$\partial_{n}V_{2}\left(  \mathcal{O}^{\prime}\right) \neq0$, but
$l\left(  h\right)  =o\left(  h^{3}\right)  $ provide the normal
derivative of an eigenfunction corresponding to $M_{2}$ vanishes
at the point $\mathcal{O}^{\prime}$. This conclusion can be
confirmed by an asymptotic analysis of eigenvalues, similar to \S
3 and \S 4. However, calculations become much more combersome and
we omit them here while refereing to \cite[Ch. 9,10]{MaNaPl} for
general asymptotic procedures.

\section{The asymptotic analysis\label{sect3}}

\subsection{The asymptotic ans\"{a}tze}

Let us examine the eigenvalues $\Lambda_{\pm}^{h}\left(
\eta\right)  $ of the spectral problem $\left(  \ref{1.14}\right)
$ in the perturbed periodicity cell $\varpi^{h}$ which are close
to the double eigenvalue $\left( \ref{1.24}\right)  $ of the
problem in $\varpi^{0}$. We fix the dual variable
of the Gel'fand transform%
\begin{equation}
\eta=\pi+\beta h^{3}, \label{2.1}%
\end{equation}
where $\beta\in\mathbb{R}$ is the deviation parameter. By varying
$\beta$, we watch over the eigenvalues $\Lambda_{\pm}^{h}\left(
\eta\right)  $ in the vicinity of the collision point in Figure
\ref{fig4},a. Note that the factor $h^{3}$ is adjusted with the
second term in the eigenvalue asymptotic
ans\"{a}tze \cite{na61} (see also \cite[Ch. 9]{MaNaPl} and \cite{NazSokAA})%
\begin{equation}
\Lambda_{\pm}^{h}\left(  \eta\right)
=\Lambda^{0}+h^{3}\Lambda_{\pm }^{^{\prime}}\left(  \beta\right)
+\widetilde{\Lambda}_{\pm}^{h}\left(
\eta\right)  . \label{2.2}%
\end{equation}
Here $\Lambda^{0}=M_{1}+\pi^{2},$ $\Lambda_{\pm}^{^{\prime}}\left(
\beta\right)  $ is a correction term to be found out and
$\widetilde{\Lambda }_{\pm}^{h}\left(  \eta\right)  $ a small
remainder to be estimated in \S \ref{sect4}. The asymptotic
ans\"{a}tze for the corresponding
eigenfunctions looks as follows:%
\begin{equation}
U_{\pm}^{h}\left(  x;\eta\right)  =U_{\pm}^{0}\left(
x;\beta\right) +h\chi\left(  x\right)  \left(  w_{\pm}^{1}\left(
\xi;\beta\right)  +hw_{\pm }^{2}\left(  \xi;\beta\right)\right)
+h^{3}U_{\pm}^{\prime}\left(  x;\beta\right)
+\widetilde{U}_{\pm}^{h}\left(  x;\eta\right). \label{2.3}%
\end{equation}
The main term%
\begin{equation}
U_{\sigma}^{0}\left(  x;\beta\right)  =a_{+}^{\sigma}\left(
\beta\right) U_{+}\left(  x\right)  +a_{-}^{\sigma}\left(
\beta\right)  U_{-}\left(
x\right)  ,\text{ }\sigma=+,-, \label{2.4}%
\end{equation}
is a linear combination of functions $\left(  \ref{1.25}\right)  $
with the coefficient column $a^{\sigma}=\left(
a_{+}^{\sigma},a_{-}^{\sigma}\right) ^{\top}$ while $\left\vert
a^{\sigma}\right\vert =1$ and $\top$ stands for transposition. The
boundary layer terms $w_{\sigma}^{q}$ are intended to compensate
for a discrepancy produced in the Dirichlet condition on the
surface $\partial\theta^{h}\cap\partial\varpi^{h}$ by the term
$U_{\pm}^{0}$. Since $w_{\sigma}\left(  \xi;\beta\right)  $ is
defined only in the set $\mathcal{U\cap\varpi}^{h}$, the cut-off
function $\chi$ is introduced in $\left(  \ref{2.3}\right)  $ such
that $\chi=0$ outside a neighborhood of $\mathcal{O}$ and $\chi=1$
in the vicinity of the point $\mathcal{O}$. The correction term
$U_{\sigma}^{\prime}$ is used to compensate for discrepancy of
$U_{\sigma}^{0}$ and $h\chi w_{\sigma}^{1}$ in the equation with
the
differential operator%
\begin{align}
\Delta_{y}+\left(  \partial_{z}+i\eta\right)
^{2}+\Lambda_{\sigma}^{h}\left(
\eta\right)   &  \sim\Delta_{y}+\left(  \partial_{z}+i\pi\right)  ^{2}%
+\Lambda^{0}+h^{3}\left(  2i\beta\left(  \partial_{z}+i\pi\right)
+\Lambda_{\sigma}^{^{\prime}}\left(  \beta\right)  \right)  +...\label{2.5}\\
&  =:L^{0}\left(  \nabla_{x}\right)  +h^{3}L^{\prime}\left(  \nabla_{x}%
;\beta\right)  +...\nonumber
\end{align}
which is decomposed in accordance with $\left(  \ref{2.1}\right)
$ and $\left(  \ref{2.2}\right)  .$ Notice that the second
boundary layer term $h\chi w_{\pm}^{2}$ is linear in
$\varkappa\left(  \mathcal{O}^{\prime }\right)  $ (cf. $\left(
\ref{1.29}\right)  $), however it does not influences
$\Lambda_{\pm}^{^{\prime}}\left(  \beta\right)  $ in $\left(
\ref{2.2}\right)  $ and becomes important only in \S \ref{sect4}
for justification estimates (see Section \ref{sect3}.4). This
observation displays the effect of opening the gap to be
independent of the curvature $\varkappa$ of the contour
$\partial\omega.$

The function $U_{\pm}^{^{\prime}}$ in $\left(  \ref{2.4}\right)  $
gets a singularity at the point $\mathcal{O}$ and, hence, we ought
to introduce another cut-off function $X_{h}$ into $\left(
\ref{2.3}\right)  $ (see $\left(  \ref{4.8}\right)  $). However,
since the asymptotic analysis in this section is formal, we avoid
to multiply $U_{\pm}^{^{\prime}}$ with $X_{h}$ here. To accept
this mathematical licence, one can assume that the coordinate
origin $\xi=0$ lies inside $\theta$, i.e.
$\mathcal{O\notin}\overline {\mathcal{\varpi}^{h}}$ for any
$h\in\left(  0,h_{0}\right]  .$

\subsection{Calculating the asymptotic terms}

Since the eigenfunction $V_{1}$ of the problem $\left(
\ref{1.4}\right)  $ is smooth near the boundary $\partial\omega$,
the Taylor formula and the
Dirichlet condition yield%
\begin{gather}
U_{\sigma}^{0}\left(  x;\beta\right)  =\left(
a_{+}^{\sigma}\left(
\beta\right)  +a_{-}^{\sigma}\left(  \beta\right)  \right)  n\partial_{n}%
V_{1}\left(  \mathcal{O}^{\prime}\right)  +%
{\displaystyle\sum\limits_{\pm}}
a_{\pm}^{\sigma}\left(  \beta\right)  \mathbf{U}^{\pm}\left(
n,s,z\right) +O\left(  \left\vert n\right\vert ^{3}+\left\vert
s\right\vert ^{3}\right)
=\label{2.6}\\
=h\left(  a_{+}^{\sigma}\left(  \beta\right)
+a_{-}^{\sigma}\left( \beta\right)  \right)
\xi_{1}\partial_{n}V\left(  \mathcal{O}^{\prime
}\right)  +h^{2}%
{\displaystyle\sum\limits_{\pm}}
a_{\pm}^{\sigma}\left(  \beta\right)  \mathbf{U}^{\pm}\left(
\xi\right)
+O\left(  h^{3}\right)  ,\text{ }\nonumber\\
\ \ \ \ \ \ \ \ \ \ \ \ \ \ \ \ \ \ \ \ \ \ \ \ \ \ \ \ \ \ \ \ \
\ \ \ \ \ \ \ \ \ \ \ \ \ \ \ \ \ \ \ \ \ \ \ \ \ \ \ \ \ \ \ \ \
\ \ \ \ \ \ \ \ \ \ \ \ \ \ \ \ \ \ \ \ \ \ \ \ \ \ \ \ \ \ \ \ \
\ \ x\in
\partial\theta^{h}\cap\partial\varpi^{h},\nonumber
\end{gather}%
\begin{equation}
\mathbf{U}^{\pm}\left(  \xi\right)  =\frac{\xi_{1}^{2}}{2}\partial_{n}%
^{2}V_{1}\left(  \mathcal{O}^{\prime}\right)
+\xi_{1}\xi_{2}\partial _{s}\partial_{n}V_{1}\left(
\mathcal{O}^{\prime}\right)  -2\pi i\delta
_{\pm,-}\xi_{1}\xi_{3}\partial_{n}V_{1}\left(
\mathcal{O}^{\prime}\right)
.\label{2.600}%
\end{equation}
Here we used the definitions of $U_{\sigma}^{0}$ and $\xi$ in
$\left( \ref{2.4}\right)  $ and $\left(  \ref{1.25}\right)  $.
Recalling the special solution $W$ of the limit problem $\left(
\ref{1.30}\right)  $ with $g\left( \xi\right)  =-\xi_{1}$, we set
\begin{equation}
w_{\sigma}^{1}\left(  \xi\right)  =\left(  a_{+}^{\sigma}\left(
\beta\right) +a_{-}^{\sigma}\left(  \beta\right)  \right)
\partial_{n}V_{1}\left(
\mathcal{O}^{\prime}\right)  W\left(  \xi\right)  \label{2.7}%
\end{equation}
in order to compensate for the main discrepancy $O\left(  h\right)
$ in $\left(  \ref{2.6}\right)  $. By the asymptotic expansion
$\left(
\ref{1.32}\right)  $ we obtain%
\begin{equation}
hw_{\sigma}^{1}\left(  \xi\right)  =-hA_{\sigma}\left(
\beta\right) \frac{\xi_{1}}{2\pi\left\vert \xi\right\vert
^{3}}+O\left(  \frac {h}{\left\vert \xi\right\vert ^{3}}\right)
=-h^{3}A_{\sigma}\left( \beta\right)  \frac{n}{2\pi r^{3}}+O\left(
\frac{h^{4}}{r^{3}}\right)
\label{2.8}%
\end{equation}
where $r=\sqrt{n^{2}+s^{2}+z^{2}}=h\left\vert \xi\right\vert $ and
\begin{equation}
A_{\sigma}\left(  \beta\right)  =\left(  a_{+}^{\sigma}\left(
\beta\right) +a_{-}^{\sigma}\left(  \beta\right)  \right)
\partial_{n}V_{1}\left(
\mathcal{O}^{\prime}\right)  P_{\theta}.\label{2.9}%
\end{equation}
After applying the differential operator $\left(  \ref{2.5}\right)
$ to the right-hand side of $\left(  \ref{2.3}\right)  $ we
collect coefficients on
$h^{3}$ and derive the differential equation%
\begin{align}
L^{0}\left(  \nabla_{x}\right)  U_{\sigma}^{^{\prime}}\left(
x;\beta\right)
& =F_{\sigma}^{\prime}\left(  x;\beta\right)  :=\label{2.10}\\
& :=-L^{\prime}\left(  \nabla_{x};\beta\right)
U_{\sigma}^{^{\prime}}\left( x;\beta\right)  +L^{0}\left(
\nabla_{x}\right)  \left(  \chi\left(  x\right)
A_{\sigma}\left(  \beta\right)  \frac{n}{2\pi r^{3}}\right)  ,\text{ }%
x\in\varpi^{0},\nonumber
\end{align}
which is to be supplied with the following Dirichlet condition on
the lateral
side $\gamma^{0}$ of the cylindrical cell $\varpi^{0}:$%
\begin{equation}
U_{\sigma}^{^{\prime}}\left(  x;\beta\right)  =0,\text{
}x\in\gamma
^{0}.\label{2.11}%
\end{equation}
The first term $L^{\prime}U_{\sigma}^{0}$ on the right of $\left(
\ref{2.10}\right)  $ is smooth in $\overline{\varpi^{0}}$ but the
second one gets a singularity at the point
$\mathcal{O}\in\gamma^{0}$. By means of $\left(  \ref{1.28}\right)
$ (see also $\left(  \ref{1.29}\right)  $ and
$\left(  \ref{2.5}\right)  $) we conclude the representation%
\begin{equation}
L^{0}\left(  \nabla_{x}\right)  =\frac{\partial^{2}}{\partial n^{2}}%
+\frac{\partial^{2}}{\partial s^{2}}+\frac{\partial^{2}}{\partial z^{2}%
}-2n\varkappa\left(  s\right)  \frac{\partial^{2}}{\partial s^{2}}%
+\mathcal{L}^{0}\left(  x,\nabla_{x}\right)  \label{2.12}%
\end{equation}
where $\mathcal{L}^{0}$ is a first-order differential operator. Hence,%
\begin{equation}
L^{0}\left(  \nabla_{x}\right)  \left(  \chi\left(  x\right)
r^{-3}n\right)
=O\left(  r^{-3}\right)  ,\text{ }r\rightarrow+0.\label{2.13}%
\end{equation}
The strong singularity $\left(  \ref{2.13}\right)  $ of the
right-hand side does not allow for a solution
$U_{\sigma}^{\prime}$ of problem $\left( \ref{2.10}\right)  ,$
$\left(  \ref{2.11}\right)  $ in the Sobolev space
$\mathring{H}^{1}\left(  \varpi^{0};\gamma^{0}\right)  .$

\subsection{The regular correction term for the eigenfunctions.}

Let us move into the scale of Kondratiev spaces
$V_{\tau}^{l}\left( \varpi^{0}\right)  $\ (see \cite{Ko} and,
e.g., \cite{NaPl, KoMaRo}) equipped
with the weighted norm%
\begin{equation}
\left\Vert U;V_{\tau}^{l}\left(  \varpi^{0}\right)  \right\Vert
=\left( \sum_{k=0}^{l}\left\Vert
\rho^{\tau-l-k}\nabla_{x}^{k}U;L^{2}\left(
\varpi^{0}\right)  \right\Vert ^{2}\right)  ^{\frac{1}{2}} \label{2.14}%
\end{equation}
where $\rho\left(  x\right)  =$dist$\left(  x,\mathcal{O}\right)
$, $\nabla_{x}^{k}U$ is the family of all order $k$ derivatives of
$U$ while
$l\in\left\{  0,1,..\right\}  $ and $\tau\in\mathbf{%
\mathbb{R}
}$ are the smoothness and weight indices, respectively. By the
one-dimensional
Hardy inequality with the particular exponent $\alpha=1,$%
\begin{equation}%
{\displaystyle\int\limits_{0}^{+\infty}}
\rho^{\alpha-1}\left\vert u\left(  \rho\right)  \right\vert
^{2}d\rho\leq
\frac{4}{\alpha^{2}}%
{\displaystyle\int\limits_{0}^{+\infty}}
\rho^{\alpha+1}\left\vert \frac{du}{d\rho}\left(  \rho\right)
\right\vert ^{2}d\rho,\text{ }\alpha>0\text{, }u\in
C_{c}^{1}\left[  0,+\infty\right)  ,
\label{2.15}%
\end{equation}
we obtain that%
\begin{equation}
\left\Vert \rho^{-1}U;L^{2}\left(  \varpi^{0}\right)  \right\Vert
^{2}\leq
c\left\Vert U;H^{1}\left(  \varpi^{0}\right)  \right\Vert ^{2}. \label{2.16}%
\end{equation}
Thus, the space%
\[
\mathring{V}_{0,per}^{1}\left(  \varpi^{0},\gamma^{0}\right)
=\left\{  U\in V_{0}^{1}\left(  \varpi^{0}\right)  :U=0\text{ on
}\gamma^{0},\text{ }U\text{ is }1-\text{periodic in }z\right\}
\]
coincides with the space $\mathring{H}_{per}^{1}\left(
\varpi^{0},\gamma
^{0}\right)  $ algebraically and topologically. This means that the mapping%
\begin{equation}
\mathring{V}_{0,per}^{1}\left(  \varpi^{0},\gamma^{0}\right)
\rightarrow \mathring{V}_{0,per}^{1}\left(
\varpi^{0},\gamma^{0}\right)  ^{\ast},
\label{2.17}%
\end{equation}
associated with the problem $\left(  \ref{2.10}\right)  $, $\left(
\ref{2.11}\right)  $, inherits all properties of the mapping%
\begin{equation}
\mathring{H}_{per}^{1}\left(  \varpi^{0},\gamma^{0}\right)
\rightarrow \mathring{H}_{per}^{1}\left(
\varpi^{0},\gamma^{0}\right)  ^{\ast}.
\label{2.18}%
\end{equation}
Moreover, theorems in \cite{Ko} on lifting smoothness and shifting
the weight indices (see also \cite[Theorems 4.1.2 and
4.2.1]{NaPl}) convey these
properties to the mapping%
\begin{equation}
\mathring{V}_{\tau,per}^{1}\left(  \varpi^{0},\gamma^{0}\right)
\cap V_{l+\tau,per}^{l+1}\left(  \varpi^{0}\right)  \rightarrow
V_{l+\tau
,per}^{l-1}\left(  \varpi^{0}\right)  \label{2.19}%
\end{equation}
in the case%
\begin{equation}
l\in\mathbf{%
\mathbb{N}
}\text{, }\tau\in\left(  -\frac{3}{2},\frac{3}{2}\right)  . \label{2.30}%
\end{equation}

\begin{remark}
\label{rem2.1} The bound $\pm\frac{3}{2}$ for the weight index
$\tau$ in $\left(  \ref{2.30}\right)  $ may be computed as
follows: the "linear"
function $x\longmapsto\chi\left(  x\right)  n$ belongs to $V_{l+\tau}%
^{l+1}\left(  \varpi^{0}\right)  $ under the restriction
$\tau>-\frac{3}{2}$ while the Poisson kernel $\chi\left(  x\right)
nr^{-3}$ lives outside $V_{l+\tau}^{l+1}\left(  \varpi^{0}\right)
$ in the case $\tau<\frac{3}{2}$. An explanation of such a
mnemonic rule, maintained by the general theory, can be found in
the introductory chapters of books \cite{NaPl, KoMaRo}$.$
\end{remark}

Taking into account the singularity of $F_{\sigma}^{\prime}$ at
the point
$\mathcal{O},$ we see that%
\[
F_{\sigma}^{\prime}\in V_{l+\tau,per}^{l-1}\left(
\varpi^{0}\right)  \text{ for any }\tau\in\left(
\frac{1}{2},\frac{3}{2}\right)  .
\]
Recall that $\Lambda^{0}=M_{1}+\pi^{2}$ is a double eigenvalue of
problem $\left(  \ref{1.14}\right)  $ in the cylindrical cell
$\varpi^{0}$ (see $\left(  \ref{1.24}\right)  $ and $\left(
\ref{1.25}\right)  $). Thus, the co-kernel of the mapping $\left(
\ref{2.19}\right)  $ is spanned over the functions $U_{\pm}$ and
the problem $\left(  \ref{2.10}\right)  $, $\left(
\ref{2.11}\right)  $ admits a solution in
$V_{l+\tau,per}^{l+1}\left( \varpi^{0}\right)  $ with
$\tau\in\left(  \frac{1}{2},\frac{3}{2}\right)  $ if
and only if%
\begin{equation}%
{\displaystyle\int\nolimits_{\varpi^{0}}}
\overline{U_{\pm}\left(  x\right)  }F_{\sigma}^{\prime}\left(
x;\beta\right)
dx=0. \label{2.31}%
\end{equation}
Note that $U_{\pm}\left(  x\right)  =O\left(  \left\vert
n\right\vert \right) $ in $\Upsilon\times\left(
-\frac{1}{2},\frac{1}{2}\right)  $ and, in view of $\left(
\ref{2.13}\right)  $, the integral in $\left(  \ref{2.31}\right)
$ is convergent.

Let the compatibility conditions $\left(  \ref{2.31}\right)  $ be
satisfied.
The orthogonality conditions%
\begin{equation}%
{\displaystyle\int\nolimits_{\varpi^{0}}}
\overline{U_{\pm}\left(  x\right)  }U_{\sigma}^{\prime}\left(
x;\beta\right)
dx=0 \label{33.0}%
\end{equation}
make the solution unique.

\begin{remark}
\label{rem3.3}. The general results \cite{Ko, MaPl2}, (see also \cite[Ch.2.3]%
{NaPl}) furnish an asymptotic form of the solution
$U_{\sigma}^{\prime}.$ By $\left(  \ref{1.28}\right)  $, $\left(
\ref{1.29}\right)  $ and $\left( \ref{2.5}\right)  $, $\left(
\ref{2.10}\right)  $, we have
\begin{equation}
F_{\sigma}^{\prime}\left(  x;\beta\right)  =r^{-5}Q_{F}\left(
n,s,z\right)
+O\left(  r^{-2}\right)  \text{, }r\rightarrow0^{+}\text{,} \label{33.1}%
\end{equation}
where $Q_{...}$ stands for a homogeneous polynomial of degree 2. A
routine and
traditional calculation brings the expansion%
\begin{equation}
U_{\sigma}^{\prime}\left(  x;\beta\right)  =r^{-3}Q_{U}\left(
n,s,z\right)
+O\left(  r^{0}\right)  \text{, }r\rightarrow0^{+}\text{,} \label{33.2}%
\end{equation}
which as well as $\left(  \ref{33.1}\right)  $ can be
differentiated under the convention $\nabla_{x}O\left(
r^{t}\right)  =O\left(  r^{t-1}\right)  $. We
need not explicit formulas for $Q_{F}$ and $Q_{U}$, however the estimate%
\begin{equation}
\nabla_{x}^{k}U_{\sigma}^{\prime}\left(  x;\beta\right)  \leq
c_{k}\left(
1+\left\vert \beta\right\vert \right)  r^{-1-k},\text{ }k=0,1..., \label{33.3}%
\end{equation}
inherited from $\left(  \ref{33.2}\right)  $, will be useful in \S \ref{sect4}%
. Constants $c_{k}$ in $\left(  \ref{33.3}\right)  $ do not depend
on the parameter $\beta$ while $Q_{F}$ and $Q_{U}$ are linear in
$\beta$ (see formula for $L^{\prime}$ in $\left(  \ref{2.5}\right)
$).
\end{remark}

All the above conclusions, of course, are known explicitly for the
Poisson kernel.

\subsection{The correction term for the eigenvalues.}

Let us compute the left-hand side of $\left(  \ref{2.31}\right)
$. Applying formulae $\left(  \ref{2.5}\right)  $, $\left(
\ref{1.000}\right)  $ and
$\left(  \ref{1.25}\right)  $, $\left(  \ref{2.4}\right)  $, we readily get%
\begin{gather}
I_{1}=%
{\displaystyle\int\nolimits_{\varpi^{0}}}
\overline{U_{\pm}\left(  x\right)  }L^{\prime}\left(
\nabla_{x};\beta\right) U_{\sigma}^{0}\left(  x;\beta\right)
dx=\ \ \ \ \ \ \ \ \ \ \ \ \ \ \ \ \ \ \ \ \ \ \ \ \ \ \ \ \ \ \ \ \ \ \ \ \ \ \ \ \ \ \ \ \label{2.32}%
\\
=%
{\displaystyle\int\nolimits_{\omega}}
\left\vert V\left(  y\right)  \right\vert ^{2}dy%
{\displaystyle\int\nolimits_{0}^{1}}
\left(  \overline{\exp\left(  \left(  -\pi\pm\pi\right)  iz\right)
}\left( a_{\sigma}^{+}\left(  \beta\right)  \left(
-2\pi\beta+\Lambda_{\sigma
}^{\prime}\left(  \beta\right)  \right)  \right)  \right.  +\nonumber\\
+\left.  \left(  a_{\sigma}^{-}\left(  \beta\right)  \left(
2\pi\beta +\Lambda_{\sigma}^{\prime}\left(  \beta\right)  \right)
\exp\left(  -2\pi
iz\right)  \right)  \right)  dz=\nonumber\\
=a_{\sigma}^{\pm}\left(  \beta\right)  \left(
\mp2\pi\beta+\Lambda_{\sigma }^{\prime}\left(  \beta\right)
\right) .\ \ \ \ \ \ \ \ \ \ \ \ \ \ \ \ \ \ \ \ \ \ \ \ \ \ \ \ \
\ \ \ \ \ \ \ \ \ \ \ \ \ \ \ \ \ \nonumber
\end{gather}
To calculate the second integral, we employ the method
\cite{MaPl1}$.$ Using the Green formula in the domain
$\varpi^{0}\backslash\mathcal{B}_{\delta}$
where $\mathcal{B}_{\delta}=\left\{  x\in\Upsilon\times\left(  -\frac{1}%
{2},\frac{1}{2}\right)  :r<\delta\right\}  $ and $\delta>0$ is small, we have%
\begin{align}
I_{2} &  =%
{\displaystyle\int\nolimits_{\varpi^{0}}}
\overline{U_{\pm}\left(  x\right)  }L^{0}\left(  \nabla_{x}\right)
\left( \chi\left(  x\right)  A_{\sigma}\dfrac{n}{4\pi
r^{3}}\right)  dx=A_{\sigma
}\underset{\delta\rightarrow0}{\lim}%
{\displaystyle\int\nolimits_{\varpi^{0}\backslash\mathcal{B}_{\delta}}}
\overline{U_{\pm}\left(  x\right)  }L^{0}\left(  \nabla_{x}\right)
\dfrac{\chi\left(  x\right)  n}{4\pi r^{3}}dx=\label{2.33}\\
&  =-A_{\sigma}\underset{\delta\rightarrow0}{\lim}%
{\displaystyle\int\nolimits_{\partial\mathcal{B}_{\delta}\cap\varpi^{0}}}
\left(  \overline{U_{\pm}\left(  x\right)
}\dfrac{\partial}{\partial \widehat{N}}\dfrac{n}{2\pi
r^{3}}-\dfrac{n}{2\pi r^{3}}\dfrac{\partial
}{\partial\widehat{N}}\overline{U_{\pm}\left(  x\right)  }\right)
ds_{x}.\nonumber
\end{align}
Here $\dfrac{\partial}{\partial\widehat{N}}=\dfrac{\partial}{\partial N}%
+i\pi\dfrac{z}{r}$ and $N$ is the interior normal on the surface
$\partial\mathcal{B}_{\delta}\cap\varpi^{0}$. Since the gradient
operator in
the curvilinear coordinates takes the form%
\[
\left(  \frac{\partial}{\partial n},\left(  1+n\varkappa\left(
s\right) \right)  ^{-1}\frac{\partial}{\partial
s},\frac{\partial}{\partial z}\right) ,
\]
we obtain%
\[
N\left(  x\right)  =\left(  r^{2}+s^{2}\left(  \left(
1+n\varkappa\left( s\right)  \right)  ^{-2}-1\right)  \right)
^{-\frac{1}{2}}\left(  n,\left( 1+n\varkappa\left(  s\right)
\right)  ^{-1}s,z\right)  ,
\]%
\[
\dfrac{\partial}{\partial N}=\left(  r^{2}+s^{2}\left(  \left(
1+n\varkappa \left(  s\right)  \right)  ^{-2}-1\right)  \right)
^{-\frac{1}{2}}\left( n\frac{\partial}{\partial n},\left(
1+n\varkappa\left(  s\right)  \right)
^{-2}s\frac{\partial}{\partial s},z\frac{\partial}{\partial
z}\right)  .
\]
Thus, computing the limit in $\left(  \ref{2.33}\right)  $, we can
make the
changes%
\[
U_{\pm}\left(  x\right)  \longmapsto n\partial_{n}V_{1}\left(  \mathcal{O}%
^{\prime}\right)  \text{, }\dfrac{\partial}{\partial\widehat{N}}%
\longmapsto\frac{\partial}{\partial
r}=\frac{n}{r}\frac{\partial}{\partial
n}+\frac{s}{r}\frac{\partial}{\partial
s}+\frac{z}{r}\frac{\partial}{\partial z}.
\]
(cf. $\left(  \ref{1.25}\right)  $, $\left(  \ref{2.6}\right)  $).
Taking the relation $mes_{2}\left(
\partial\mathcal{B}_{\delta}\cap\varpi^{0}\right)
=2\pi\delta^{2}+O\left(  \delta^{3}\right)  $ into account, we
then arrive at
the formula%
\begin{align}
I_{2} &  =-A_{\sigma}\left(  \beta\right)  \partial_{n}V_{1}\left(
\mathcal{O}^{\prime}\right)  \underset{\delta\rightarrow0}{\lim}%
{\displaystyle\int\nolimits_{\partial\mathcal{B}_{\delta}\cap\varpi^{0}}}
\left(  n\dfrac{\partial}{\partial r}\dfrac{n}{2\pi
r^{3}}-\dfrac{n}{2\pi
r^{3}}\dfrac{\partial n}{\partial r}\right)  ds_{x}=\label{2.34}\\
&  =A_{\sigma}\left(  \beta\right)  \partial_{n}V_{1}\left(  \mathcal{O}%
^{\prime}\right)  =\left(  a_{\sigma}^{+}\left(  \beta\right)  +a_{\sigma}%
^{-}\left(  \beta\right)  \right)  P_{\theta}\left\vert \partial_{n}%
V_{1}\left(  \mathcal{O}^{\prime}\right)  \right\vert
^{2}.\nonumber
\end{align}
Here we have used notation $\left(  \ref{2.9}\right)  $ and
further we set
$\mathcal{P=}P_{\theta}\left\vert \partial_{n}V_{1}\left(  \mathcal{O}%
^{\prime}\right)  \right\vert ^{2}$ as in $\left(
\ref{1.35}\right)  .$

By $\left(  \ref{2.32}\right)  $ and $\left(  \ref{2.34}\right)
$, the compatibility conditions $\left(  \ref{2.31}\right)  $
reduce to the system of
two algebraic equations%
\[
\pm2\pi\beta a_{\sigma}^{\pm}\left(  \beta\right)
+\mathcal{P}\left( a_{\sigma}^{+}\left(  \beta\right)
+a_{\sigma}^{-}\left(  \beta\right) \right)
=\Lambda_{\sigma}^{\prime}\left(  \beta\right)  a_{\sigma}^{\pm
}\left(  \beta\right)  .
\]
Eigenvalues of the corresponding matrix%
\begin{equation}
\left(
\begin{tabular}
[c]{ll}%
$2\pi\beta+\mathcal{P}$ & $\mathcal{P}$\\
$\mathcal{P}$ & $-2\pi\beta+\mathcal{P}$%
\end{tabular}
\ \right)  \label{2.35}%
\end{equation}
look as follows%
\begin{equation}
\Lambda_{\pm}^{\prime}\left(  \beta\right)  =\mathcal{P\pm}\sqrt
{\mathcal{P}^{2}+4\pi^{2}\beta^{2}}. \label{2.36}%
\end{equation}

\subsection{The second term in the boundary layer.}

Even in the case $\varkappa\left(  \mathcal{O}^{\prime}\right)
=0$, e.g., the contour $\partial\omega$ is flat near the point
$\mathcal{O}^{\prime}$ and, by
formulas $\left(  \ref{1.1}\right)  $ and $\left(  \ref{1.28}\right)  $,%
\[
\partial_{n}^{2}V_{1}\left(  \mathcal{O}^{\prime}\right)  =-\varkappa\left(
\mathcal{O}^{\prime}\right)  \partial_{n}V_{1}\left(
\mathcal{O}^{\prime }\right)  =0,
\]
the second term $\left(  \ref{2.600}\right)  $ of the discrepancy
$\left( \ref{2.6}\right)  $ does not vanish. The Sobolev norm of
the functions $h^{q}\chi w_{\pm}^{q}$ is $O\left(
h^{p+\frac{1}{2}}\right)  $ and, therefore, our aim to derive
estimates with the bound $ch^{\frac{7}{2}}$ forces us to deal with
$h^{2}\chi w_{\pm}^{2}$ in \S 4, although this term, owing to the
proper decay as $\left\vert \xi\right\vert \rightarrow\infty$,
does not influence the correction term
$h^{3}\Lambda_{\pm}^{\prime}\left( \beta\right)  $ in $\left(
\ref{2.2}\right)  .$

The Taylor formula $\left(  \ref{2.6}\right)  $ gives immediately
the boundary
condition%
\begin{equation}
w_{\pm}^{2}\left(  \xi;\beta\right)  =-a_{\pm}^{+}\left(
\beta\right) \mathbf{U}^{+}\left(  \xi\right)  -a_{\pm}^{-}\left(
\beta\right)
\mathbf{U}^{-}\left(  \xi\right)  ,\text{ }\xi\in\partial\Theta. \label{2.40}%
\end{equation}
To derive the differential equation%
\begin{equation}
-\Delta_{\xi}w_{\pm}^{2}\left(  \xi;\beta\right)  =-\left(
\varkappa\left( \mathcal{O}^{\prime}\right)
\frac{\partial}{\partial\xi_{1}}-2\varkappa \left(
\mathcal{O}^{\prime}\right)  \xi_{1}\frac{\partial^{2}}{\partial
\xi_{2}^{2}}+2\pi i\frac{\partial}{\partial\xi_{3}}\right)  \widetilde{w}%
_{\pm}^{1}\left(  \xi;\beta\right)  ,\text{ }\xi\in\Theta, \label{2.41}%
\end{equation}
requires much more elaborated analysis based on the procedure
\cite[\S 2.2, Ch.4]{MaNaPl} of discrepancies rearrangement. First,
the differential operator $\mathcal{L}_{1}\left(
\xi,\nabla_{\xi}\right)  $ on the right of $\left(
\ref{2.41}\right)  $ comes from the expansions $\left(
\ref{1.29}\right)  $, $\left(  \ref{2.5}\right)  $ and
\[
\left(  \frac{\partial}{\partial z}+i\pi\right)  ^{2}=\frac{1}{h^{2}}%
\frac{\partial^{2}}{\partial z^{2}}+\frac{1}{h}2\pi
i\frac{\partial}{\partial z}-\pi^{2},
\]
in other words, $\mathcal{L}_{1}\left(  \xi,\nabla_{\xi}\right)  $
appears as a coefficient on $h^{-1}$ in the decomposition of
$\Delta_{y}+\left(
\partial_{z}+i\pi\right)  ^{2}$ in the stretched curvilinear coordinates $\xi
$. Second,%
\begin{equation}
\widetilde{w}_{\pm}^{1}\left(  \xi;\beta\right)
=w_{\pm}^{1}\left(  \xi ;\beta\right)  +A_{\pm}\left(
\beta\right)  \left(  2\pi\left\vert \xi\right\vert ^{3}\right)
^{-1}\xi_{1}=O\left(  \left\vert \xi\right\vert
^{-3}\right)  \label{2.414}%
\end{equation}
while, according to the rearrangement procedure mentioned above,
the main asymptotic term in $\left(  \ref{2.8}\right)  $ is
detached from the
right-hand side of $\left(  \ref{2.41}\right)  $ because the expression%
\begin{equation}
A_{\pm}\left(  \beta\right)  \mathcal{L}_{1}\left(
\xi,\nabla_{\xi}\right) \frac{\xi_{1}}{2\pi\left\vert
\xi\right\vert ^{3}}=h^{3}A_{\pm}\left(
\beta\right)  \mathcal{L}_{1}\left(  n,s,z,\partial_{n},\partial_{s}%
,\partial_{z}\right)  \frac{n}{2\pi r^{3}} \label{2.42}%
\end{equation}
has too slow decay $O\left(  \left\vert \xi\right\vert
^{-3}\right)  $ at infinity and, hence, putting $\left(
\ref{2.42}\right)  $ into $\left( \ref{2.41}\right)  $ would lead
to the insufficient decay rate of the boundary layer term. The
transmission of certain unsuitable constituents from one limit
problem to the other limit problem and the preservation of the
behavior of lower order asymptotic terms as $r\rightarrow0^{+}$
and $\left\vert \xi\right\vert \rightarrow+\infty$ implies the
absence of the rearrangement procedure. Recall that, indeed, the
expression $\left(  \ref{2.42}\right)  $ with the cut-off function
$\chi$ is a part of the right-hand side in $\left(
\ref{2.10}\right)  $, and notice that the detachment made in
$\left( \ref{2.41}\right)  $ helps crucially to derive necessary
estimates in \S \ref{sect4}.

Similarly to Remark \ref{rem3.3} one, based on a general result in
\cite{Ko} (see also \cite[\S 3.5, \S 6.4]{NaPl}), may conclude the
existence of a unique decaying solution of the problem $\left(
\ref{2.41}\right)  ,$ $\left(
\ref{2.40}\right)  $ and the relation%
\begin{equation}
w_{\pm}^{2}\left(  \xi;\beta\right)  =O\left(  \left\vert
\xi\right\vert ^{-2}\right)  \text{, }\left\vert \xi\right\vert
\rightarrow\infty.
\label{2.43}%
\end{equation}
We emphasize that, due to the factor $h^{2}$ instead of $h$ in
$w_{\pm}^{1}$, the decay rate $\left(  \ref{2.43}\right)  $ damps
down the influence of $w_{\pm}^{2}$ on $U_{\pm}^{\prime}$ (cf. a
calculation in $\left( \ref{2.8}\right)  $).

\subsection{Simple eigenvalues}

If $\eta\in\left[  0,2\pi\right)  $ and $\eta\neq\pi$, the
eigenvalues $\Lambda_{\pm}^{0}\left(  \eta\right)  =M_{1}+\left(
\eta-\pi\pm\pi\right) ^{2}$ (see $\left(  \ref{1.22}\right)  $)
are simple and the corresponding eigenfunctions $U_{\pm}^{0}$ are
still given by $\left(  \ref{1.25}\right)  .$ The asymptotic
structures remain the same as for $\left(  \ref{2.1}\right)  $
but loose dependence on the parameter $\beta$. In particular,%
\begin{equation}
w_{\pm}\left(  \xi\right)  =\partial_{n}V\left(
\mathcal{O}^{\prime}\right)
W\left(  \xi\right)  \label{2.37}%
\end{equation}
(cf. $\left(  \ref{2.7}\right)  $ with
$a_{\pm}^{\sigma}=\delta_{\sigma,\pm}$)
and $U_{\pm}^{\prime}$ satisfies the equation%
\begin{gather}
-\left(  \Delta_{y}+\left(  \partial_{z}+i\eta\right)  \right)
^{2}U_{\pm }^{\prime}\left(  x\right)  -\Lambda_{\pm}^{0}\left(
\eta\right)  U_{\pm
}^{\prime}\left(  x\right)  =\Lambda_{\pm}^{0}\left(  \eta\right)  U_{\pm}%
^{0}\left(  x\right)
+\ \ \ \ \ \ \ \ \ \ \ \ \ \ \ \ \ \ \ \ \ \ \ \ \ \ \ \ \ \ \ \ \ \ \ \ \ \ \ \ \ \ \ \ \ \ \ \ \ \ \ \label{2.38}%
\\
\ \ \ \ \ \ \ \ \ \ \ \ \ \ \ \ \ \ \ \ \ \ \ \ \ \ \ \ \ \ \ \ \
+\left( \Delta_{y}+\left(  \partial_{z}+i\eta\right)  \right)
^{2}\left(  \chi\left( x\right)  \dfrac{n}{2\pi r^{3}}\right)  ,\
\ x\in\varpi^{0},
\end{gather}
supplied with the Dirichlet conditions $\left(  \ref{2.11}\right)
$ and the periodicity conditions. Since the eigenvalue
$\Lambda_{\pm}^{0}\left( \eta\right)  $ is simple, only one
compatibility condition must be verified, and repeating the
calculation $\left(  \ref{2.33}\right)  $, $\left(
\ref{2.34}\right)  $ brings the equalities%
\begin{equation}
\Lambda_{\pm}^{\prime}\left(  \eta\right)
=\mathcal{P\colon=}P_{\theta }\left\vert \partial_{n}V_{1}\left(
\mathcal{O}^{\prime}\right)  \right\vert
^{2}.\label{2.39}%
\end{equation}
One readily sees that formulas $\left(  \ref{2.2}\right)  $,
$\left( \ref{1.22}\right)  $, $\left(  \ref{2.39}\right)  $ with
$\eta=\pi+\beta h^{3}$ bring about an expansion for the
eigenvalues $\Lambda_{\pm}^{h}\left( \pi+\beta h^{3}\right)  $
which differs from the expansion obtained in the previous section.
This lack of coincidence originates in ignoring the second
compatibility condition, namely the norm of the inverse operator
in $\varpi^{h}$ restricted onto a subspace of co-dimension $1$
grows when $\eta\rightarrow\pi$ and the eigenvalues
$\Lambda_{\pm}^{h}\left( \eta\right)  $ approach one the other.

In \S \ref{sect4} the most attention is paid for an appropriate
estimate of the remainder $\widetilde{\Lambda}_{\pm}^{h}\left(
\pi+\beta h^{3}\right)  $ in a sufficiently wide range of the
deviation parameter $\beta.$

\section{Justification of the asymptotic expansion\label{sect4}}

\subsection{The operator formulation of the cell problem.}

To estimate the asymptotic remaiders
$\widetilde{\Lambda}_{\pm}^{h}\left( \eta\right)  $ in formulas
$\left(  \ref{1.33}\right)  $ and $\left( \ref{1.30}\right)  $, we
employ the following fact which is known as "Lemma on almost
eigenvalues and eigenvectors" and can be found in, e. g.,
\cite{ViLu, BiSo} with much more general formulation.

\begin{lemma}
\label{lem4.1} Let $\mathbf{H}$ be an Hilbert space and
$\mathbf{K}$ be a
compact self-adjoint positive operator in $\mathbf{H}$. If $\mathbf{y}%
\in\mathbf{H}$ and $\varphi\in\mathbf{%
\mathbb{R}
}_{+}$ meet the conditions%
\begin{equation}
\left\Vert \mathbf{y};\mathbf{H}\right\Vert =1\text{, }\left\Vert
\mathbf{Ky}-\varphi\mathbf{y};\mathbf{H}\right\Vert
=\delta\in\left(
0,\varphi\right)  , \label{4.1}%
\end{equation}
then the segment $\left[  \varphi-\delta,\varphi+\delta\right]
\subset
\mathbf{%
\mathbb{R}
}_{+}$ contains an eigenvalue $\psi$ of the operator $\mathbf{K}$.
\end{lemma}

The space $\mathring{H}_{per}^{1}\left(
\varpi^{h};\gamma^{h}\right)  $ equipped with the scalar product
\begin{equation}
\left\langle u,v\right\rangle _{\eta}=\left(  \left(
\nabla_{x}+i\eta e_{3}\right)  U,\left(  \nabla_{x}+i\eta
e_{3}\right)  V\right)  _{\varpi^{h}}
\label{4.2}%
\end{equation}
is denoted by $\mathcal{H}\left(  \eta\right)  $. Here
$\eta\in\left[ 0,2\pi\right)  $ while the Friedrichs inequality
(cf. the middle part of $\left(  \ref{1.39}\right)  $) provides
the positiveness of the Hermitian form $\left(  \ref{4.2}\right)
$.

By a simple argument, the operator $\mathcal{K}\left(  \eta\right)
$,
determined by the identity%
\begin{equation}
\left\langle \mathcal{K}\left(  \eta\right)  U,V\right\rangle
_{\eta}=\left( U,V\right)  _{\varpi^{h}},\text{ }U\text{,
}V\in\mathcal{H}\left(
\eta\right)  , \label{4.3}%
\end{equation}
is compact, self-adjoint and positive. Owing to \cite[Thm.
9.2.1]{BiSo}, the spectrum of this operator consists of the
essential spectrum $\left\{
0\right\}  $ and the discrete spectrum%
\begin{equation}
\psi_{1}^{h}\left(  \eta\right)  \geq\psi_{2}^{h}\left(
\eta\right) \geq...\geq\psi_{p}^{h}\left(  \eta\right)
\geq...\rightarrow0^{+}.
\label{4.4}%
\end{equation}
Comparing $\left(  \ref{4.2}\right)  $, $\left(  \ref{4.3}\right)
$ with
$\left(  \ref{1.36}\right)  $, we observe the relationship%
\begin{equation}
\Lambda_{p}^{h}\left(  \eta\right)  =\psi_{1}^{h}\left(
\eta\right)  ^{-1}
\label{4.5}%
\end{equation}
between entries in the eigenvalue sequences $\left(
\ref{1.37}\right)  $ and $\left(  \ref{4.5}\right)  $.

\subsection{Approximation solutions for the spectral problem.}

Let us consider the most interesting case $\left(
\ref{2.1}\right)  $. We
shorten the notation as follows:%
\[
\mathcal{H}_{\beta}=\mathcal{H}\left(  \pi\pm\beta h^{3}\right)
,\text{ }\mathcal{K}_{\beta}=\mathcal{K}\left(  \pi\pm\beta
h^{3}\right)  ,\text{ }\left\langle \text{ },\text{ }\right\rangle
_{\pi+\beta h^{3}}=\left\langle \text{ },\text{ }\right\rangle
_{\left(  \beta\right)  }.
\]
Furthermore, we set
\begin{equation}
\varphi_{\pm}=l_{\pm}^{-1}\text{, }\mathcal{Y}_{\pm}=\left\Vert
Y_{\pm
};\mathcal{H}_{\beta}\right\Vert ^{-1}Y_{\pm}, \label{4.6}%
\end{equation}
where%
\begin{equation}
l_{\pm}=\Lambda_{0}+h^{3}\Lambda_{\pm}^{\prime}\left(
\beta\right)  ,
\label{4.7}%
\end{equation}%
\begin{align}
Y_{\pm}\left(  x\right)   &  =X_{h}\left(  x\right)
U_{\pm}^{0}\left( x;\beta\right)  +\left(  1-X_{h}\left(  x\right)
\right)  \left( n\partial_{n}U_{\pm}^{0}\left(
\mathcal{O};\beta\right)  +\mathbf{U}_{\pm
}\left(  x;\beta\right)  \right)  +\label{4.8}\\
&  +h\chi\left(  x\right)  \left(  w_{\pm}^{1}\left(
\xi;\beta\right) +hw_{\pm}^{2}\left(  \xi;\beta\right)
+h^{3}X_{h}\left(  x\right)  U_{\pm }^{\prime}\left(
x;\beta\right)  \right)  .\nonumber
\end{align}
Notice that the dependence on $h$ and $\beta$ is not indicated in
$\left( \ref{4.6}\right)  $. Addenda on the right of $\left(
\ref{4.7}\right)  $ have been determined in $\left(
\ref{2.2}\right)  $, $\left(  \ref{2.36}\right) $. However, the
function $\left(  \ref{4.8}\right)  $ has still to be specified.
First, the regular terms $U_{\pm}^{0},$ $U_{\pm}^{\prime}$ and the
boundary layer $w_{\pm}$ were constructed in \S \ref{sect2} while
the
coefficient column $a^{\pm}$ in the linear combination $\left(  \ref{2.4}%
\right)  $ had to be an eigenvector of the matrix $\left(
\ref{2.35}\right)
,$%
\begin{align}
a^{\pm}  &  =a_{0}^{-\frac{1}{2}}\left(  \mathcal{P},-2\pi\beta\pm
\sqrt{\mathcal{P}^{2}+4\pi^{2}\beta^{2}}\right)  ,\label{4.9}\\
a_{0}  &  =2\left(  \mathcal{P}^{2}+4\pi^{2}\beta^{2}\mp2\pi\beta
\sqrt{\mathcal{P}^{2}+4\pi^{2}\beta^{2}}\right)  ,\nonumber\\
\left(  a^{\pm},a^{\pm}\right)  _{\mathbf{%
\mathbb{R}
}^{2}}  &  =1,\left(  a^{\pm},a^{\mp}\right)  _{\mathbf{%
\mathbb{R}
}^{2}}=0.\nonumber
\end{align}
Second, the cut-off function $X_{h}$ is determined in $\left(  \ref{1.41}%
\right)  $ while, according to (\ref{2.600}), we set%
\begin{equation}
\mathbf{U}_{\pm}\left(  x;\beta\right)  =a_{\pm}^{+}\left(
\beta\right) \mathbf{U}^{+}\left(  n,s,z\right)
+a_{\pm}^{-}\left(  \beta\right)
\mathbf{U}^{-}\left(  n,s,z\right)  . \label{4.444}%
\end{equation}
Third, $X_{h}$ cuts off the regular terms near the cavern
$\theta^{h}$ and,
thus, due to the relation (see section \ref{sect3}.3)%
\[
w_{\pm}^{1}\left(  \xi;\beta\right)
=-\xi_{1}\partial_{n}U_{\pm}^{0}\left( \mathcal{O}\right)
,\xi\in\partial\Theta,
\]
and the boundary condition (\ref{2.40}) for $w_{\pm}^{2}\left(
\xi ;\beta\right)  $, the function $Y_{\pm}$ vanishes on the
surface $\gamma^{h}$
and, therefore, falls into $\mathring{H}_{per}^{1}\left(  \varpi^{h}%
;\gamma^{h}\right)  .$ We finally mention that $X_{h}$ smooths
down the correction term $U_{\pm}^{\prime}$ which gets a
singularity at $\mathcal{O}$ (see Section \ref{sect3}.3).

Calculating the norm $\left\Vert
\mathbf{U}_{\pm};\mathcal{H}_{\beta }\right\Vert $, we obtain

\begin{align}
\left\langle U_{\pm}^{0},U_{\pm}^{0}\right\rangle _{\pi+\beta h^{3}%
}=&\left\Vert \left(  \nabla_{x}+i\left(  \pi+\beta h^{3}\right)
\right)
U_{\pm}^{0};L^{2}\left(  \varpi^{h}\right)  \right\Vert ^{2}=\nonumber\\
=&\left\Vert \left(  \nabla_{x}%
+i\pi\right)  U_{\pm}^{0};L^{2}\left(  \varpi^{0}\right)
\right\Vert
^{2}+R_{\pm}^{0}=\nonumber\\
 =&\left\Vert
\nabla_{y}V;L^{2}\left( \omega\right)  \right\Vert
^{2}+\pi^{2}\left\Vert V;L^{2}\left(
\omega\right)  \right\Vert ^{2}+R_{\pm}^{0}=\nonumber\\
 =&M_{1}+\pi^{2}+R_{\pm}^{0},\nonumber\\
\left\vert R_{\pm}^{0}\right\vert \leq & cr_{\pm}%
^{0}\left(  1+r_{\pm}^{0}\right),\label{4.10}\\
 r_{\pm}^{0}=&\left\Vert \nabla_{x}U_{\pm}%
^{0};L^{2}\left(  \varpi\backslash\varpi^{h}\right)  \right\Vert
+\left\vert \pi+\beta h^{3}\right\vert \left\Vert
U_{\pm}^{0};L^{2}\left(  \varpi
\backslash\varpi^{h}\right)  \right\Vert +\nonumber\\
&+\left\vert
\beta\right\vert h^{3}\left\Vert U_{\pm}^{0};L^{2}\left(
\varpi^{h}\right)  \right\Vert \leq ch^{\frac{3}{2}}\left(
1+\left\vert \beta\right\vert h^{\frac{3}{2}}\right) \nonumber
\end{align}

and%
\begin{align}
\Vert \left(  1-X_{h}\right) & \left(  U_{\pm}^{0}\left(
\cdot ,\beta\right)  -n\partial_{n}U_{\pm}^{0}\left(
\mathcal{O};\beta\right)
-\mathbf{U}_{\pm}\right)  ;\mathcal{H}_{\beta}\Vert \leq\label{4.11}\\
 &\leq c\left(
{\displaystyle\int\nolimits_{0}^{2c_{X}h}}
\left(  r^{4}+\left(  h^{-2}+\left\vert \pi+\beta h^{3}\right\vert
^{2}\right)  r^{6}\right)  r^{2}dr\right)  ^{\frac{1}{2}}\leq\nonumber\\
&\leq ch^{\frac{7}{2}}\left(  1+\left\vert \beta\right\vert
h^{5}\right) ,\nonumber
\end{align}%
\begin{align}
h^{p}\left\Vert \chi w_{\pm}^{p};\mathcal{H}_{\beta}\right\Vert  &
\leq ch^{p}\left(  h^{\frac{1}{2}}\left\Vert
\nabla_{\xi}w_{\pm};L^{2}\left( \Theta\right)  \right\Vert +\left(
1+\left\vert \beta\right\vert h^{3}\right)
h^{\frac{3}{2}}\left\Vert w_{\pm};L^{2}\left(  \Theta\right)
\right\Vert \right)  \leq\label{4.12}\\
&  \leq ch^{p+\frac{1}{2}}\left(  1+\left\vert \beta\right\vert
h^{4}\right) ,\ \ \ \ \ \ p=1,2,\nonumber
\end{align}%
\begin{align}
h^{3}\left\Vert
X_{h}U_{\pm}^{\prime};\mathcal{H}_{\beta}\right\Vert  & \leq
ch^{3}\left(
{\displaystyle\int\nolimits_{c_{X}h}^{\text{diam}\left(
\varpi^{0}\right)  }}
\left(  r^{-4}+\left(  h^{-2}+\left\vert \pi\beta P\beta
h^{3}\right\vert
^{2}\right)  r^{-2}\right)  r^{2}dr\right)  ^{\frac{1}{2}}\leq\label{4.13}\\
& \leq ch^{\frac{7}{2}}\left(  1+\left\vert \beta\right\vert
h^{3}\right) .\nonumber
\end{align}
Let us comment on the above calculations. In $\left(
\ref{4.10}\right)  $ and then in $\left(  \ref{4.11}\right)  $ we
applied the explicit formulas $\left(  \ref{2.4}\right)  ,$
$\left(  \ref{4.9}\right)  ,$ $\left( \ref{1.000}\right)  $ and
also the relations $\left(  \ref{2.6}\right)  ,$ $\left(
\ref{1.41}\right)  .$ The inequalities $\left(  \ref{4.12}\right)
$ hold true due to the coordinate dilation $x\mapsto\xi$ and the
inclusion $w_{\pm}^{p}\in H^{1}\left(  \Theta\right)  $ inherited
from the expansion $\left(  \ref{1.32}\right)  $ and the relation
$\left(  \ref{2.43}\right)  $. Finally, the term
$h^{3}U_{\pm}^{\prime}$ with $X_{h}$ was treated by means of the
estimates $\left(  \ref{33.3}\right)  $ taking the properties
$\left( \ref{1.41}\right)  $ of the cut-off function into account.

Imposing the restriction
\begin{equation}
\left\vert \beta\right\vert \leq h^{-\frac{5}{4}}\beta_{0}, \label{4.15}%
\end{equation}
which damps down the parameter $\beta$ in all bounds in the
inequalities $\left(  \ref{4.10}\right)  -\left(
\ref{4.12}\right)  $. We emphasize that the weaker restriction
$\left\vert \beta\right\vert \leq\beta_{1}h^{-\frac
{3}{2}}$ is sufficient here (see the last estimate in $\left(  \ref{4.10}%
\right)  $), however in the sequel we need $\left(
\ref{4.15}\right)  $ and just this restriction has been imposed in
Theorem \ref{teo1}. We observe $\left(  \ref{4.8}\right)  -\left(
\ref{4.13}\right)  $ and conclude that, for a small $\beta_{0}>0$
in $\left(  \ref{4.15}\right)  ,$ the following
inequality is valid:%
\begin{equation}
\left\Vert Y_{\pm};\mathcal{H}_{\beta}\right\Vert
\geq\frac{1}{2}\left(
M_{1}+\pi^{2}\right)  . \label{4.16}%
\end{equation}
Moreover, by $\left(  \ref{4.7}\right)  $ and $\left(
\ref{2.36}\right)  ,$ under the same condition $\left(
\ref{4.15}\right)  ,$ the numbers
(\ref{4.7}) are subject to%
\begin{equation}
l_{\pm}\geq\frac{1}{2}\left(  M_{1}+\pi^{2}\right)  . \label{4.17}%
\end{equation}

\subsection{Justifying the asymptotic expansions of eigenvalues.}

For the approximate solution $\left(  \ref{4.6}\right)  $, the
quantity
$\delta$ in $\left(  \ref{4.1}\right)  $ takes the form%
\begin{align}
\delta &  =\left\Vert
\mathcal{K}_{\beta}\mathcal{Y}_{\pm}-\varphi_{\pm
}\mathcal{Y}_{\pm};\mathcal{H}_{\beta}\right\Vert =\sup\left\vert
\left\langle
\mathcal{K}_{\beta}\mathcal{Y}_{\pm}-\varphi_{\pm}\mathcal{Y}_{\pm
},V\right\rangle _{\left(  \beta\right)  }\right\vert =\label{4.21}\\
&  =l_{\pm}^{-1}\left\Vert Y_{\pm};\mathcal{H}_{\beta}\right\Vert ^{-1}%
\sup\left\vert \left\langle
Y_{\pm}-l_{\pm}\mathcal{K}_{\beta}Y_{\pm
},V\right\rangle _{\left(  \beta\right)  }\right\vert =\nonumber\\
&  =l_{\pm}^{-1}\left\Vert Y_{\pm};\mathcal{H}_{\beta}\right\Vert ^{-1}%
\sup\left\vert S_{\pm}\left(  V\right)  \right\vert \leq
c\sup\left\vert S_{\pm}\left(  V\right)  \right\vert ,\ \nonumber
\end{align}
where the supremum is calculated over all functions
$V\in\mathcal{H}_{\beta}$ such that $\left\Vert
V;\mathcal{H}_{\beta}\right\Vert =1$ and
\begin{align}
S_{\pm}\left(  V\right)    & =\left(  \left(  \nabla_{x}+i\left(
\eta+\beta h^{3}\right)  e_{3}\right)  Y_{\pm},\left(
\nabla_{x}+i\left(  \pi+\beta
h^{3}\right)  \right)  \nabla_{x}V\right)  _{\varpi^{h}}-\label{4.22}\\
& -\left(  \Lambda_{0}+h^{3}\Lambda_{\pm}^{\prime}\left(
\beta\right) \right)  \left(  Y_{\pm},V\right)
_{\varpi^{h}}.\nonumber
\end{align}
Notice that the Friedrichs and Hardy inequalities (see $\left(  \ref{1.39}%
\right)  $ and $\left(  \ref{2.15}\right)  $ with $\alpha=1$)
provide the
estimate%
\begin{equation}
\left\Vert \nabla_{x}V;L^{2}\left(  \varpi^{h}\right)  \right\Vert
+\left\Vert \rho^{-1}V;L^{2}\left(  \varpi^{h}\right)  \right\Vert
\leq c\left\Vert
V;\mathcal{H}_{\beta}\right\Vert =c.\label{4.23}%
\end{equation}
We extend the test function $V$ by null onto $\varpi$ and subtract
from
$\left(  \ref{4.22}\right)  $ the following scalar products:%
\begin{equation}
\left(  \left(  \nabla_{x}+i\pi e_{3}\right)  U_{\pm}^{0},\left(
\nabla
_{x}+i\pi e_{3}\right)  V\right)  _{\varpi}-\Lambda^{0}\left(  U_{\pm}%
^{0},V\right)  _{\varpi}=0,\label{4.24}%
\end{equation}%
\begin{equation}
h^{2}\left(
\nabla_{\xi}\widetilde{w}_{\pm}^{1},\nabla_{\xi}\left(  \chi
V\right)  \right)  _{\varpi}=0,\label{4.25}%
\end{equation}%
\begin{equation}
h^{3}\left(  \nabla_{\xi}w_{\pm}^{2},\nabla_{\xi}\left(  \chi
V\right) \right)  _{\Theta}-h^{3}\left(
\mathcal{L}_{1}\widetilde{w}_{\pm}^{1},\chi
V\right)  _{\Theta}=0,\label{4.255}%
\end{equation}%
\begin{equation}%
\begin{array}
[c]{l}%
\left(  \left(  \nabla_{x}+i\pi e_{3}\right)
U_{\pm}^{\prime},\left( \nabla_{x}+i\pi e_{3}\right)  \left(
X_{h}V\right)  \right)  _{\varpi
}-\Lambda^{0}\left(  U_{\pm}^{\prime},X_{h}V\right)  _{\varpi}-\\
\\
\ \ \ \ \ \ \ \ \ \ \ \ \ \ \ \ \ -\Lambda_{\pm}^{\prime}\left(
\beta\right) \left(  U_{\pm}^{0},X_{h}V\right)
_{\varpi}+2i\beta\left(  U_{\pm}^{\prime },\left(
\partial_{z}+i\pi\right)  \left(  X_{h}V\right)  \right)  _{\varpi
}-\\
\\
\ \ \ \ \ \ \ \ \ \ \ \ \ \ \ \ \ \ -\Lambda_{\pm}^{\prime}\left(
\beta\right)  \left(  \left(  \nabla_{x}+i\pi e_{3}\right)
\chi\left(  2\pi r^{3}\right)  ^{-1}n,\left(  \nabla_{x}+i\pi
e_{3}\right)  \left( X_{h}V\right)  \right)  _{\varpi}=0.
\end{array}
\label{4.26}%
\end{equation}
The equality $\left(  \ref{4.24}\right)  $ is just the integral
identity $\left(  \ref{1.14}\right)  $. The function $\chi V$ is
written in the stretched curvilinear coordinates $\xi$ (see
$\left(  \ref{1.26}\right)  $), vanishes on $\partial\Theta$ and
has a compact support; thus $\left( \ref{4.255}\right)  $ and
$\left(  \ref{4.25}\right)  $ follow from the equation $\left(
\ref{2.41}\right)  $ and the harmonicity of the function
$\left(  \ref{2.414}\right)  $, respectively. Finally, $\left(  \ref{4.26}%
\right)  $ is but a consequence of $\left(  \ref{2.10}\right)  $,
$\left( \ref{2.11}\right)  $; note that the test function $X_{h}V$
vanishes near the point $\mathcal{O}$ where $U_{\pm}^{\prime}$ has
the strong singularity $\left(  \ref{33.2}\right)  $.

In the next two sections we estimate terms which are left in
$\left( \ref{4.22}\right)  $ after subtracting left-hand sides of
$\left( \ref{4.24}\right)  -\left(  \ref{4.26}\right)  $ and
obtain the common bound $ch^{\frac{7}{2}}\left\Vert
V;\mathcal{H}_{\beta}\right\Vert $. By virtue of $\left(
\ref{4.16}\right)  ,$ $\left(  \ref{4.17}\right)  $ and $\left(
\ref{4.21}\right)  $, Lemma \ref{lem4.1} delivers an eigenvalue
$\psi_{\pm
}^{h}$ of the operator $\mathcal{K}_{\beta}$ such that%
\[
\left\vert \psi_{\pm}^{h}-l_{\pm}\right\vert \leq
c_{\psi}h^{\frac{7}{2}}.
\]
Using $\left(  \ref{4.5}\right)  $ and $\left(  \ref{4.7}\right)
$, this
formula yields%
\begin{align}
\left\vert \Lambda_{\pm}^{h}\left(  \pi+\beta h^{3}\right)  -\Lambda^{0}%
-h^{3}\Lambda_{\pm}^{\prime}\left(  \beta\right)  \right\vert  &
\leq c_{\psi}h^{\frac{7}{2}}\Lambda_{\pm}^{h}\left(  \pi+\beta
h^{3}\right) \left(  \Lambda^{0}+h^{3}\Lambda_{\pm}^{\prime}\left(
\beta\right)  \right)
,\label{4.27}\\
\Lambda_{\pm}^{h}\left(  \pi+\beta h^{3}\right)  \left(
1-c_{\psi}h^{\frac {7}{2}}\left(
\Lambda^{0}+h^{3}\Lambda_{\pm}^{\prime}\left(  \beta\right)
\right)  \right)   &
\leq\Lambda^{0}+h^{3}\Lambda_{\pm}^{\prime}\left( \beta\right)
.\nonumber
\end{align}
Thus, recalling the condition $\left(  \ref{4.15}\right)  $ and
choosing $h_{0}>0$ such that the factor $\Lambda_{\pm}^{h}\left(
\pi+\beta h^{3}\right)  $ on the left of $\left(
\ref{4.27}\right)  $ is bigger than
$\dfrac{1}{2}$ we arrive at the estimate%
\[
\left\vert \Lambda_{\pm}^{h}\left(  \pi+\beta h^{3}\right)  -\Lambda^{0}%
-h^{3}\Lambda_{\pm}^{\prime}\left(  \beta\right)  \right\vert \leq
c_{\Lambda }h^{\frac{7}{2}},
\]
which proves Theorem \ref{teo1}.

\subsection{Discrepancies of the regular terms.}

Proceeding with $U_{\pm}^{0}$, we have to take into account the scalar product%
\[
I_{1}=h^{6}\beta^{2}\left(  U_{\pm}^{0},V\right)  _{\varpi^{0}}%
\]
and the cut-off function $X_{h}$ in $\left(  \ref{4.8}\right)  $ resulting in%
\[
I_{2}=\left(  \left(  \nabla_{x}+i\left(  \pi+\beta h^{3}\right)
e_{3}\right)  \left(  1-X_{h}\right)  \left(
U_{\pm}^{0}-n\partial_{n}U_{\pm }^{0}-\mathbf{U}_{\pm}\right)
,\left(  \nabla_{x}+i\left(  \pi+\beta h^{3}\right)  e_{3}\right)
V\right)  _{\varpi^{0}}.
\]
Other constituents of $\left(  \ref{4.22}\right)  $, involving
$U_{\pm}^{0},$ are included to either $\left(  \ref{4.24}\right)
$, or $\left( \ref{4.26}\right)  $. Recall that the test function
$V$ is extended by zero on the whole cell $\varpi^{0}$. In view of
formulas $\left(  \ref{2.3}\right)  $, $\left(  \ref{4.9}\right)
$ and $\left(  \ref{4.23}\right)  $, a bound for
$\left\vert I_{1}\right\vert $ looks as follows:%
\[
ch^{6}\beta^{2}\left\Vert V;L^{2}\left(  \varpi^{0}\right)
\right\Vert \leq
ch^{\frac{7}{2}}\left(  h^{\frac{5}{4}}\beta\right)  ^{2}\leq ch^{\frac{7}{2}%
}.
\]
Here we have used the restriction $\left(  \ref{4.15}\right)  $ on
the parameter $\beta$, which also applies in further calculations.
In the sequel we skip mentioning this argument.

By $\left(  \ref{2.6}\right)  $, $\left(  \ref{1.41}\right)  $ and
$\left(
\ref{4.23}\right)  $, we have%
\[
\left\vert I_{2}\right\vert \leq c\left(
{\displaystyle\int\nolimits_{0}^{2c_{X}h}}
\left(  r^{4}+\left(  h^{-2}\left(  1+\left\vert \beta\right\vert
h^{3}\right)  ^{2}\right)  r^{6}\right)  r^{2}dr\right)  ^{\frac{1}{2}%
}\left\Vert V;\mathcal{H}_{\beta}\right\Vert \leq
ch^{\frac{7}{2}}.
\]
We now consider the terms due to the transportation of $X_{h}$
from $U_{\pm
}^{\prime}$ to $V$, namely%
\[
I_{3}=h^{3}\left(  U_{\pm}^{\prime}\nabla_{x}X_{h},\left(
\nabla_{x}+i\left( \pi+\beta h^{3}\right)  e_{3}\right)  V\right)
_{\varpi^{0}}+h^{3}\left( \left(  \nabla_{x}+i\left(  \pi+\beta
h^{3}\right)  e_{3}\right)  U_{\pm
}^{\prime},V\nabla_{x}X_{h}\right)  _{\varpi^{0}}.\ \ \ \
\]
We obtain%
\begin{gather*}
\left\vert I_{3}\right\vert \leq ch^{3}\left(  h^{-2}%
{\displaystyle\int\nolimits_{c_{x}h}^{2c_{x}h}}
r^{-2}r^{2}dr\left\Vert V;\mathcal{H}_{\beta}\right\Vert
^{2}+\right.
\ \ \ \ \ \ \ \ \ \ \ \ \ \ \ \ \ \ \ \ \ \ \ \ \ \ \ \ \ \ \ \ \ \ \ \ \ \ \ \ \ \ \ \ \ \ \ \ \ \ \ \ \ \ \ \ \ \ \ \ \ \ \ \ \ \ \ \\
\left.  +%
{\displaystyle\int\nolimits_{c_{x}h}^{2c_{x}h}}
\left(  r^{-4}+\left(  1+\left\vert \beta\right\vert h^{3}\right)  ^{2}%
r^{-2}\right)  r^{2}drh^{-2}\left\Vert \left(  h^{-1}r\right)  ^{-1}%
V;L^{2}\left(  \varpi^{0}\right)  \right\Vert ^{2}\right)  ^{\frac{1}{2}}%
\leq\\
\leq ch^{3}\left(  h^{-2}h+h^{-1}h^{2}h^{-2}\right)  \left\Vert V;\mathcal{H}%
_{\beta}\right\Vert ^{2}\leq ch^{\frac{7}{2}}%
.\ \ \ \ \ \ \ \ \ \ \ \ \ \ \ \ \ \ \ \ \ \ \ \ \ \ \ \ \ \ \ \ \
\ \ \ \ \ \ \ \ \ \ \ \ \ \ \ \ \ \ \ \ \
\end{gather*}
Here we used the estimates $\left(  \ref{33.3}\right)  $ and
$\left( \ref{4.23}\right)  $ for $U_{\pm}^{\prime}$ and $V$,
respectively, while the factors $h^{-2}$ and $h^{2}$ are caused by
the differentiation of the cut-off function and the relation
$c_{X}\leq h^{-1}r\leq2c_{X}$ on supp$\left\vert
\nabla_{x}X_{h}\right\vert $ (see $\left(  \ref{1.41}\right)  $
and compare
with $\left(  \ref{1.444}\right)  $). The list of other remaining terms reads%
\[
i\beta h^{6}\left(  \left(  U_{\pm}^{\prime},\left(
\nabla_{x}+i\left( \pi+\beta h^{3}\right)  e_{3}\right)
X_{h}V\right)  _{\varpi^{0}}-\left( \left(  \nabla_{x}+i\pi
e_{3}\right)  U_{\pm}^{\prime},X_{h}V\right) _{\varpi^{0}}\right)
\]%
\[
-h^{6}\Lambda_{\pm}^{\prime}\left(  \beta\right)  \left(
U_{\pm}^{\prime },X_{h}V\right)  _{\varpi^{0}}.
\]
Estimates for these terms with the bound $ch^{\frac{7}{2}}$ become
evident
after applying the inequality%
\[
\left\Vert X_{h}V;\mathcal{H}_{\beta}\right\Vert \leq c\left\Vert
V;\mathcal{H}_{\beta}\right\Vert
\]
following from $\left(  \ref{1.41}\right)  $ and $\left(
\ref{4.23}\right)  $.

\subsection{Discrepancies of the boundary layer terms.}

First of all, we replace $w_{\pm}^{1}$ by
$\widetilde{w}_{\pm}^{1}$ since the main asymptotic term,
subtracted in $\left(  \ref{2.414}\right)  $ from the boundary
layer solution $w_{\pm}^{1}$, has been included into the equation
$\left(  \ref{2.10}\right)  $ and, therefore, the expression
$\left( \ref{4.26}\right)  $.

Next, the inequality%
\begin{gather*}
h\left\vert \left(  \left(
\widetilde{w}_{\pm}^{1}+hw_{\pm}^{2}\right) \nabla_{x}\chi,\left(
\nabla_{x}+i\left(  \pi+\beta h^{3}\right) e_{3}\right)  V\right)
_{\varpi^{0}}\right.
+\ \ \ \ \ \ \ \ \ \ \ \ \ \ \ \ \ \ \ \ \ \ \ \ \ \ \ \ \ \ \ \ \ \ \ \ \ \ \ \ \ \ \ \ \ \ \ \ \ \ \ \ \ \ \ \ \ \ \ \ \ \ \ \ \ \ \ \ \ \ \ \ \ \ \ \ \ \ \ \ \\
+\left.  \left(  \left(  \nabla_{x}+i\left(  \pi+\beta
h^{3}\right) e_{3}\right)  \left(
\widetilde{w}_{\pm}^{1}+hw_{\pm}^{2}\right)
,V\nabla_{x}\chi\right)  _{\varpi^{0}}\right\vert \leq\\
\leq ch\left(
{\displaystyle\int\nolimits_{\text{supp}\left\vert
\nabla_{x}\chi\right\vert }}
\left(  \left\vert \frac{r}{h}\right\vert ^{6}+h^{2}\left\vert \frac{r}%
{h}\right\vert ^{-4}\right)  dx\right)  ^{\frac{1}{2}}\left\Vert
V;\mathcal{H}_{\beta}\right\Vert \leq ch^{4}%
\ \ \ \ \ \ \ \ \ \ \ \ \ \ \ \ \ \ \ \ \ \ \ \ \ \ \ \ \ \ \ \ \
\ \ \ \ \ \ \ \ \ \ \ \ \ \
\end{gather*}
permits for transporting the cut-off function from
$\widetilde{w}_{\pm}^{p}$ to $V$. Note that derivatives of $\chi$
vanish in a neighborhood of the point
$\mathcal{O}$ and the integral over the support of $\left\vert \nabla_{x}%
\chi\right\vert $ converges. Moreover, the relations $\left(  \ref{2.414}%
\right)  $ and $\left(  \ref{2.43}\right)  $ for
$\widetilde{w}_{\pm}^{1}$ and $w_{\pm}^{2}$, respectively, can be
differentiated in the case $h^{-1}\xi\in $supp$\left\vert
\nabla_{x}\chi\right\vert $ that was used in the inequality.

Finally, we write%
\begin{gather}
h\left\vert \left(  i\beta h^{3}\widetilde{w}_{\pm}^{1},\left(
\partial
_{z}+i\left(  \pi+\beta h^{3}\right)  \right)  \chi V\right)  _{\varpi^{0}%
}\right\vert +h\left\vert \left(  \left(  \partial_{z}+i\left(
\pi+\beta h^{3}\right)  \right)  \widetilde{w}_{\pm}^{1},i\beta
h^{3}\chi V\right)
_{\varpi^{0}}\right\vert +\label{4.27b}\\
\ \ \ \ \ \ \ \ \ \ \ \ \ \ \ +h^{2}\left\vert \left(  \left(
i\left( \pi+\beta h^{3}\right)  w_{\pm}^{2},\left(
\partial_{z}+i\left(  \pi+\beta h^{3}\right)  \right)  \chi
V\right)  \right)  _{\varpi^{0}}\right\vert
+\nonumber\\
\ \ \ \ \ \ \ \ \ \ \ \ \ \ \ +h^{2}\left\vert \left(  \left(
\partial _{z}+i\left(  \pi+\beta h^{3}\right)  \right)
w_{\pm}^{2},i\left(  \pi+\beta
h^{3}\right)  \chi V\right)  _{\varpi^{0}}\right\vert \leq\nonumber\\
\leq ch^{2}\left(  1+\left\vert \beta\right\vert h^{2}\right)
\left( h^{\frac{3}{2}}\left(  \left\Vert
\widetilde{w}_{\pm}^{1};L^{2}\left( \Theta\right)  \right\Vert
+\left\Vert w_{\pm}^{2};L^{2}\left(  \Theta\right) \right\Vert
^{2}\left\Vert V;\mathcal{H}_{\beta}\right\Vert \right)  +\right.
\nonumber\\
\ \ \ \ \ \ \ \ \ \ \ \ \ \ \ \ \ \ \ \ \ \ \left.
+h^{\frac{3}{2}}\left(
\left\Vert \left\vert \xi\right\vert \nabla_{\xi}\widetilde{w}_{\pm}^{1}%
;L^{2}\left(  \Theta\right)  \right\Vert +\left\Vert \left\vert
\xi\right\vert \nabla_{\xi}w_{\pm}^{2};L^{2}\left(  \Theta\right)
\right\Vert \right) \left\Vert r^{-1}V;L^{2}\left(
\varpi^{0}\right)  \right\Vert \right)
\leq\nonumber\\
\leq ch^{\frac{7}{2}}%
,\ \ \ \ \ \ \ \ \ \ \ \ \ \ \ \ \ \ \ \ \ \ \ \ \ \ \ \ \ \ \ \ \
\ \ \ \ \ \ \ \ \ \ \ \ \ \ \ \ \ \ \ \ \ \ \ \ \ \ \ \ \ \ \ \ \
\ \ \ \ \ \ \ \ \ \ \ \ \ \ \ \ \ \ \ \ \ \ \ \ \ \ \ \ \
\nonumber
\end{gather}%
\begin{align*}
&  h\left\vert \left(
\Lambda_{0}+h^{3}\Lambda_{\pm}^{\prime}\left( \beta\right)
\right)  \left(  \widetilde{w}_{\pm}^{1}+hw_{\pm}^{2},\chi
V\right)  _{\varpi^{0}}\right\vert \leq\\
&  \leq ch\left(  h^{\frac{5}{2}}\left\Vert \left\vert
\xi\right\vert \widetilde{w}_{\pm}^{1};L^{2}\left(  \Theta\right)
\right\Vert \left\Vert
r^{-1}V;L^{2}\left(  \varpi^{0}\right)  \right\Vert +hh^{\frac{3}{2}%
}\left\Vert w_{\pm}^{2};L^{2}\left(  \Theta\right)  \right\Vert
\left\Vert V;L^{2}\left(  \varpi^{0}\right)  \right\Vert \right)
\leq ch^{\frac{7}{2}}.
\end{align*}
Here we have made the transform $x\mapsto\xi$ which brings the
factor $h^{\frac{3}{2}}$ on the $L^{2}-$norms of
$\widetilde{w}_{\pm}^{1}$, $\left\vert \xi\right\vert
\nabla_{\xi}\widetilde{w}_{\pm}^{1}$, $\widetilde{w}_{\pm}^{2},$
$\left\vert \xi\right\vert \nabla_{\xi}w_{\pm}^{2}$ and the factor
$h^{\frac{5}{2}}$ on $\left\vert \xi\right\vert \widetilde
{w}_{\pm}^{1}$. We emphasize that all norms in $L^{2}\left(
\Theta\right)  $ figuring in $\left(  \ref{4.27b}\right)  $ appear
to be finite due to the relations $\left(  \ref{2.414}\right)  $
and $\left(  \ref{2.43}\right)  $.

The above considerations demonstrate that the inner products
involving boundary layer components in $\left(  \ref{4.22}\right)
$, can be changed with the error $O\left(  h^{\frac{7}{2}}\right)
$ for the sum of the
following integrals%
\begin{equation}
i\pi h%
{\displaystyle\int\nolimits_{0}^{1}}
{\displaystyle\int\nolimits_{\partial\omega}}
{\displaystyle\int\nolimits_{0}^{d}}
\left(  \widetilde{w}_{\pm}^{1}\overline{\partial_{z}\left(  \chi
V\right) }-\partial_{z}\widetilde{w}_{\pm}^{1}\overline{\chi
V}\right)  \left(
1+n\varkappa\right)  dndsdz,\label{4.28}%
\end{equation}%
\begin{equation}
h%
{\displaystyle\int\nolimits_{0}^{1}}
{\displaystyle\int\nolimits_{\partial\omega}}
{\displaystyle\int\nolimits_{0}^{d}}
\left(
\partial_{n}\widetilde{w}_{\pm}^{1}\overline{\partial_{n}\left(
\chi V\right)  }+\left(  1+n\varkappa\right)
^{-2}\partial_{s}\widetilde{w}_{\pm
}^{1}\overline{\partial_{s}\left(  \chi V\right)
}+\partial_{z}\widetilde {w}_{\pm}^{1}\overline{\partial_{z}\left(
\chi V\right)  }\right)  \left(
1+n\varkappa\right)  dndsdz,\label{4.29}%
\end{equation}%
\begin{equation}
h^{2}%
{\displaystyle\int\nolimits_{0}^{1}}
{\displaystyle\int\nolimits_{\partial\omega}}
{\displaystyle\int\nolimits_{0}^{d}}
\left(  \partial_{n}w_{\pm}^{2}\overline{\partial_{n}\left(  \chi
V\right) }+\left(  1+n\varkappa\right)
^{-2}\partial_{s}w_{\pm}^{2}\overline {\partial_{s}\left(  \chi
V\right)  }+\partial_{z}w_{\pm}^{2}\overline {\partial_{z}\left(
\chi V\right)  }\right)  \left(  1+n\varkappa\right)
dndsdz,\ \ \ \ \ \ \ \ \ \ \ \ \ \ \ \ \label{4.30}%
\end{equation}
where $1+n\varkappa\left(  s\right)  $ is the Jacobian, the
differential
operator $\nabla_{x}$ in the curvilinear coordinates takes the form%
\[
\left(  \frac{\partial}{\partial n},\left(  1+n\varkappa\left(
s\right) \right)  ^{-1}\frac{\partial}{\partial
s},\frac{\partial}{\partial z}\right)
\]
and $d>0$ is chosen such that supp$\chi$ belongs to the $d-$
neighborhood of the contour $\partial\omega$ allowing for the
curvilinear coordinate system $\left(  n,s,z\right)  $.

Replacing $1+n\varkappa\left(  s\right)  $ by $1$ in $\left(  \ref{4.30}%
\right)  $ brings an error which does not exceed%
\[
ch^{2}\left\vert
{\displaystyle\int\nolimits_{\varpi^{0}}}
\left\vert \nabla_{x}\left(  \chi\left(  x\right)  V\left(
x\right)  \right) \right\vert n\left\vert
h^{-1}\nabla_{\xi}w_{\pm}^{2}\left(  \xi\right) \right\vert
dx\right\vert \leq ch^{2}h^{\frac{3}{2}}\left\Vert \left\vert
\xi\right\vert \nabla_{\xi}w_{\pm}^{2}\left(  \xi\right)
;L^{2}\left( \Theta\right)  \right\Vert \left\Vert
V;\mathcal{H}_{\beta}\right\Vert \leq ch^{\frac{7}{2}}.\ \
\]
The resultant integral (with $n\varkappa=0$ in $\left(
\ref{4.30}\right)  $) turns into the first addendum in $\left(
\ref{4.255}\right)  $ by the transform $\left(  n,s,z\right)
\longmapsto\xi.$

The same procedure works for $\left(  \ref{4.28}\right)  $ with an
error less
than%
\begin{align*}
&  ch%
{\displaystyle\int\nolimits_{\varpi^{0}}}
n\left(  \left\vert \widetilde{w}_{\pm}^{1}\left(  \xi\right)
\right\vert \left\vert \nabla_{x}\left(  \chi\left(  x\right)
V\left(  x\right)  \right) \right\vert +r\left\vert
h^{-1}\nabla_{\xi}\widetilde{w}_{\pm}^{1}\left( \xi\right)
\right\vert r^{-1}\left\vert \chi\left(  x\right)  V\left(
x\right)  \right\vert \right)  dx\\
&  \leq ch\left(  h^{\frac{5}{2}}\left\Vert \left\vert
\xi\right\vert \widetilde{w}_{\pm}^{1}\left(  \xi\right)
;L^{2}\left(  \Theta\right)
\right\Vert +h^{\frac{5}{2}}\left\Vert \left\vert \xi\right\vert ^{2}%
\nabla_{\xi}\widetilde{w}_{\pm}^{1}\left(  \xi\right)
;L^{2}\left( \Theta\right)  \right\Vert \right)  \left\Vert
V;\mathcal{H}_{\beta }\right\Vert \leq ch^{\frac{7}{2}}\ \
\end{align*}
while the norms in $L^{2}\left(  \Theta\right)  $ still stay finite in view of%
\[
\left\vert \xi\right\vert \left\vert \widetilde{w}_{\pm}^{1}\left(
\xi\right)  \right\vert =O\left(  \left\vert \xi\right\vert
^{-2}\right) \text{, }\left\vert \xi\right\vert ^{2}\left\vert
\nabla_{\xi}\widetilde {w}_{\pm}^{1}\left(  \xi\right)
\right\vert =O\left(  \left\vert \xi\right\vert ^{-2}\right)
\]
(cf. $\left(  \ref{2.414}\right)  $). The resultant integral becomes%
\begin{equation}
\pi ih^{3}%
{\displaystyle\int\nolimits_{\Theta}}
\left(  \widetilde{w}_{\pm}^{1}\frac{\partial\overline{\chi
V}}{\partial \xi_{3}}-\overline{\chi
V}\frac{\partial\widetilde{w}_{\pm}^{1}}{\partial
\xi_{3}}\right)  d\xi=-2\pi ih^{3}%
{\displaystyle\int\nolimits_{\Theta}}
\overline{\chi
V}\frac{\partial\widetilde{w}_{\pm}^{1}}{\partial\xi_{3}}d\xi.
\label{4.31}%
\end{equation}
In $\left(  \ref{4.29}\right)  $ we substitute
$1+h\xi_{1}\varkappa\left(
\mathcal{O}^{\prime}\right)  $ and $1-h\xi_{1}\varkappa\left(  \mathcal{O}%
^{\prime}\right)  $ for $1+n\varkappa\left(  s\right)  $ and
$\left( 1+n\varkappa\left(  s\right)  \right)  ^{-1}$,
respectively. The concomitant error again gets order
$h^{\frac{7}{2}}$ but, in addition to the expression on
the left of $\left(  \ref{4.25}\right)  $, we obtain the integral%
\begin{align}
&  h^{3}\varkappa\left(  \mathcal{O}^{\prime}\right)
{\displaystyle\int\nolimits_{\Theta}}
\xi_{1}\left(  \dfrac{\partial\widetilde{w}_{\pm}^{1}}{\partial\xi_{1}}%
\dfrac{\partial\overline{\chi V}}{\partial\xi_{1}}-\dfrac{\partial
\widetilde{w}_{\pm}^{1}}{\partial\xi_{2}}\dfrac{\partial\overline{\chi V}%
}{\partial\xi_{2}}+\dfrac{\partial\widetilde{w}_{\pm}^{1}}{\partial\xi_{3}%
}\dfrac{\partial\overline{\chi V}}{\partial\xi_{3}}\right)  d\xi\label{4.32}\\
&  =-h^{3}\varkappa\left(  \mathcal{O}^{\prime}\right)
{\displaystyle\int\nolimits_{\Theta}}
\left(  \Delta_{\xi}\widetilde{w}_{\pm}^{1}+\xi_{1}\dfrac{\partial
\widetilde{w}_{\pm}^{1}}{\partial\xi_{1}}-2\dfrac{\partial^{2}\widetilde
{w}_{\pm}^{1}}{\partial\xi_{2}^{2}}\right)  \chi Vd\xi.\nonumber
\end{align}
Since $\widetilde{w}_{\pm}^{1}$ is a harmonics, the integrals
$\left( \ref{4.32}\right)  $ and $\left(  \ref{4.31}\right)  $,
according to $\left( \ref{2.41}\right)  $, form the second term on
the left of (\ref{4.255}).

We have verified the fact which had been announced in the end of
Section \ref{sect4}.3. Our proof of Theorem \ref{teo1} is now
completed.

\bigskip

\textbf{Acknowledgements.} This paper was prepared during the
visit of S.A. Nazarov to Department of Civil Engineering of Second
University of Naples and it was supported by project "Asymptotic
analysis of composite materials and thin and non-homogeneous
structures" (Regione Campania, law n.5/2005) and by the grant RFFI
- 06--01--257.

\end{document}